\documentclass[12pt]{amsart} 

\usepackage{amssymb} 
\usepackage{enumerate, amsfonts, latexsym}
  
\newtheorem{theorem}{Theorem}[section]

\newtheorem{lemma} [theorem]{Lemma}
\newtheorem{proposition}[theorem]{Proposition}
\newtheorem{corollary} [theorem] {Corollary}

\newtheorem{observation}[theorem] {Observation}

\newtheorem{convention}[theorem]{Convention}
\newtheorem{workingassumption}[theorem]{Working Assumption}

\theoremstyle{definition}
\newtheorem{definition}[theorem]{Definition}
\newtheorem{example}[theorem]{Example}

\theoremstyle{remark}
\newtheorem{remark}[theorem]{Remark}
\newtheorem{remarks}[theorem]{Remarks}

\numberwithin{equation}{section}

\def\<{\langle}
\def\>{\rangle}
\def\|{{\ |\ }}

\def\-{\underline}

\def\e{\varepsilon} 
 
\def\N{\mathbb N}

\def\On{\text{\rm{Out}} (F_n)}
\def\Out{${\rm{Out}} (F_n)$}
\def\Aut{${\rm{Aut}} (F_n)$}

\def\gep{{\rm{GEP}}}
\def\pep{{\rm{$\Psi$EP}}}
\def\atom{atom}

\def\ssm{\smallsetminus}

 \def\NF {Nibbled Future}
 \def\nf {nibbled future}
 
 \def\nib {nibbling}

\catcode`\@=11
\def\serieslogo@{\relax}
\def\@setcopyright{\relax}
\catcode`\@=12

\begin{document}

\title[beaded decompositions]
{Free-group automorphisms, train tracks, and the
beaded decomposition}

\author[Martin R. Bridson]{Martin R.~Bridson}
\address{Martin R.~Bridson\\
Mathematics Department\\
180 Queen's Gate\\
London, SW7 2BZ\\
U.K. }
\email{m.bridson@ic.ac.uk} 

\author[Daniel Groves]{Daniel Groves}
\address{Daniel Groves\\
Department of Mathematics\\
California Institute of Technology\\
Pasadena, CA 9110, USA }
\email{groves@caltech.edu}

\date{27 July 2006.}


\thanks{The first author's work was supported in part
by Research Fellowships from the EPSRC and a Royal Society
Wolfson Research Merit Award.  The second author was supported
in part by a Junior Research Fellowship at Merton College, Oxford, and by NSF Grant DMS-0504251. 
We thank these organisations for their support.}

\begin{abstract} {We study the automorphisms $\phi$ of a finitely generated
free group $F$. Building on the train-track
technology of Bestvina, Feighn and Handel, we provide a topological
representative $f:G\to G$ of a power of $\phi$ that  behaves very much
like the realization on the rose of a positive automorphism. This
resemblance is encapsulated in the {\em Beaded Decomposition Theorem}
which describes the structure of paths in $G$ obtained by repeatedly
passing to $f$-images of an edge and taking subpaths. 
This decomposition is the key 
to adapting our proof of the quadratic isoperimetric inequality for 
$F\rtimes_\phi\mathbb Z$, with $\phi$ positive, to the general case.}
\end{abstract}

\maketitle

The study of automorphisms of free groups is  informed greatly by the  
analogies with automorphisms of free-abelian groups and surface groups, but one
often has to work considerably harder in the free group case in order to
obtain the appropriate analogues of familiar results from these other
contexts. Nowhere is this more true than in the quest for suitable normal
forms and geometric representatives. One can gain insight into the nature of
individual elements of $\hbox{GL}(n,\mathbb Z)$ by realizing them
as diffeomorphisms  of the $n$-torus. Likewise, one analyzes individual
elements   of the mapping class group by realizing them
as diffeomorphisms  of a surface.  The situation for \Aut\ and \Out\ 
is more complicated: the natural choices of classifying space $K(F_n,1)$ 
are finite graphs of genus $n$, and no element of infinite order in \Out\ 
is induced by the action on $\pi_1(Y)$ of a 
homeomorphism of $Y$. Thus the best that
one can hope for in this situation is to identify a graph $Y_\phi$ that admits a 
homotopy equivalence inducing $\phi$ and that has additional structure
well-adapted to $\phi$. 
This is the context of the {\em train track technology} of
 Bestvina, Feighn and Handel \cite{BH2, BFH, BFH2}.

Their work results in a decomposition
theory for elements of \Out\ that is closely analogous to (but 
more complicated than) the Nielsen-Thurston theory for surface
automorphisms.  The finer features of 
the topological normal forms that they obtain are adapted to the problems that
they wished to solve in each of their papers: the Scott conjecture in \cite{BH2}
and the Tits alternative in the series of papers \cite{BFH, BFH2, BFHGeomDed}. 
The problem that we wish to solve in this series of papers (of which \cite{BGroves} is the first, this is the second and \cite{bg3} is the third), is that of determining
the Dehn functions of  all free-by-cyclic groups. This requires
a further refinement of the train-track technology. Specifically,
we must adapt our topological representatives 
so as to make tractable the problem of determining the isoperimetric properties
of the mapping torus of the homotopy equivalence $f: Y_\phi \to Y_\phi$ realizing
an iterate of $\phi$. 

The prototype for a train-track representative is the obvious realization
of a positive automorphism on the rose. This motivates the following
strategy for the Dehn-function problem: first we proved the theorem
in the case of positive automorphisms \cite{BGroves}, where one already
encounters most of the web of large-scale cancellation phenomena that explain
why the general theorem is true; then, in the general case, we
follow the architecture of the proof from  \cite{BGroves}, using a suitably refined
train-track description of the automorphism in place of the positivity assumption.
We shall
see in \cite{bg3} that this approach works remarkably well. However, in
order to bring it to fruition,  
 one must deal with myriad additional complexities 
arising from the intricacies of cancellation that do not arise in the positive case.  

The  properties of the
topological representative $f:G\to G$ constructed in  \cite{BFH}  allow one
to control the manner in which a path $\sigma$ evolves as
one looks at its iterated images under $f$, and one might naively suppose that this
is the key issue that one must overcome in translating the arguments from the
positive case \cite{BGroves} to the general case \cite{bg3}. However, upon
close inspection one discovers this is actually only a fraction  of the story, the
point being that when a corridor evolves in the time flow on a van Kampen diagram,
the interaction of the forward iterates of the individual edges is such that the
basic {\em splitting} of paths established in \cite{BFH} gets broken. It is to overcome
this difficulty that we need the notion of {\em hard splitting}
introduced in Section \ref{s:hard}; such splittings are denoted $\sigma_1\odot\sigma_2$.

In the analysis of van Kampen diagrams, the class of ``broken" paths that one must
understand are the residues of the images of a single edge that survive
repeated cancellation during the corridor flow.
 In the language of the topological representative
$f:G\to G$, this amounts to understanding {\em monochromatic
paths}, as defined below. 
Every edge-path $\rho$ in $G$ admits a unique maximal splitting 
into edge paths;
our purpose in this article is to understand the nature of factors in the
case where $\rho$ is monochromatic (grouping certain of the factors into
larger units).

To this end, we identify a small number of basic 
units into which the iterated images of 
monochromatic paths split; the key feature of this splitting is that it 
is robust enough   to withstand the
difficulties caused by cancellation in van Kampen diagrams. The basic units
are defined so as to ensure that they enjoy those features of
individual edges that proved important in the positive case \cite{BGroves}. 
We call the units {\em beads}. The vocabulary of beads is as follows.

Let $f:G\to G$ be a topological representative and let $f_\#(\sigma)$
denote the tightening rel endpoints of the image of an edge-path $\sigma$.
Following \cite{BH2}, if $f_\#(\tau)=\tau$  we call $\tau$ a
{\em Nielsen} path.
A path $\rho$ in $G$ is called a {\em growing exceptional path} (\gep) if
either $\rho$ or $\bar{\rho}$ is of
the form $E_i\bar{\tau}^k \bar{E_j}$ where $\tau$ is a Nielsen
path, $k \geq 1$, $E_i$ and $E_j$ are parabolic edges,  $f(E_i) = E_i \odot \tau^m$, 
$f(E_j) = E_j \odot\tau^n$, and $n > m > 0$. If it is $\rho$ (resp.~$\overline\rho$) 
that is of this form, then proper initial (resp.~terminal)
sub edge-paths of $\rho$ are called \pep s
({\em pseudo-exceptional paths}). 

Let $f: G \to G$ be an improved relative train track map and $r, J \ge 1$ integers.
Then {\em $r$-monochromatic} paths in $G$ are defined by a simple  recursion:
edges in $G$ are $r$-monochromatic and if $\rho$ is  an
$r$-monochromatic path then every sub edge-path of
$f^r_\#(\rho)$ is $r$-monochromatic.
A {\em $(J,f)$-atom} is an 
$f$-monochromatic edge-path of length at most $J$ that admits no non-vacuous
hard splitting into edge-paths.
 
An edge-path $\rho$ is {\em $(J,f)$-beaded} if it admits a hard splitting $\rho = \rho_1 \odot \cdots \odot \rho_k$ where each $\rho_i$ is a \gep, a \pep, a $(J,f)$-atom, or
an indivisible Nielsen path of length at most $J$ (where \gep s, \pep s and Nielsen paths are defined with respect to the map $f$).

The following is the main result of this paper. 

\medskip{\bf{Beaded Decomposition Theorem:}} {\em For every $\phi\in{\rm{Out}}(F_n)$,
there exist positive integers $k, r$ and $J$ such that $\phi^{k}$ has an improved
relative train-track representative $f:G\to G$ with the property that every $r$-monochromatic path in $G$ is $(J,f)$-beaded.}
\medskip 

(See Subsection \ref{ss:Iterate} below for a precise description of what we mean by
the map $f_{\#}^r$.)
 
As is clear from the preceding   discussion, our main motivation for
developing the Beaded  Decomposition is its application in \cite{bg3}. 
The import of the current paper in \cite{bg3} has been deliberately distilled
into  this single statement, and the technical addenda in Section \ref{refinements},
so that a reader who is willing to accept these
as articles of faith can proceed directly from \cite{BGroves} to \cite{bg3}.

We expect that our particular refinement of the train-track technology may prove
useful in other contexts.  This expectation stems from
the general point that the development of refined topological representatives
leads to insights into purely algebraic questions about free-group
automorphisms.  See \cite{BG-Growth} for a concrete illustration of this.\footnote{\cite{BG-Growth} contains some results about the growth of words under iterated
automorphisms.  A previous version of this paper contained an incorrect version of these results.  We thank Gilbert Levitt for pointing out our error.}

\setcounter{tocdepth}{1}
\tableofcontents

\section{Improved relative train track maps} \label{TrainTracks}

In this section we collect and refine those elements of the train-track
technology that we shall need. Most of the material here is drawn
directly from  \cite{BH2} and \cite{BFH}. 

The philosophy behind
train tracks is to find an {\em efficient} topological representative for
an outer automorphism of $F$.  Precisely what it means for a graph map
to be {\em efficient} is spelled out in this section.

\subsection{Paths, Splittings, Turns and Strata}
Let $G$ be a graph. Following \cite{BFH}, we try to
reserve the term {\em path} for a
map $\sigma : [0,1] \to G$ that is either constant or an immersion (i.e.
{\em tight}).  The reverse path $t\mapsto \sigma(1-t)$ will be denoted $\overline \sigma$.
We
conflate the map $\sigma$ with its monotone reparameterisations (and
even its image, when this does not cause confusion).  Given an arbitrary continuous map
$\rho : [0,1] \to G$, we denote by $[\rho]$ the unique (tight)
{\em path}  homotopic
rel  endpoints to $\rho$. In keeping with the notation of the previous
section, given $f:G\to G$ and a path $\sigma$  in $G$,  we write $f_{\#}(\sigma)$
to denote $[f(\sigma)]$. We are primarily concerned
with {\em edge paths},  i.e. those paths $\sigma$ for which $\sigma(0)$
and $\sigma(1)$ are vertices.

We consider only maps $f:G\to G$ that send vertices to vertices and  
edges to edge-paths (not necessarily to single edges). If there is 
an isomorphism  $F \cong \pi_1 G$ such that $f$ induces $\mathcal O
\in \text{\rm{Out}}(F)$, then one says that $f$ {\em represents} $\mathcal O$.

\noindent{\bf Notation:} Given a map $f$ of graphs and a path $\rho$
in the domain, we'll follow the  standard practice of denoting by
$f_\#(\rho)$  the unique
locally-injective edge path that is homotopic rel endpoints to  $f(\rho)$. 

\subsection{Replacing $f$ by an Iterate} \label{ss:Iterate}
In order to obtain good topological representatives of outer
automorphisms, one has to replace the given map by a large
iterate. It is important to be clear what one means by {\em iterate}
in this context, since   we wish to consider only
topological  representatives  whose restriction to
each edge is an immersion and this property is not inherited by
(naive) powers of the map.

Thus we deem the phrase\footnote{and obvious
variations on it}
{\em replacing $f$ by an iterate}, to mean that for fixed
$k\in\mathbb N$, we pass from consideration of $f : G \to G$
to consideration of  the map $f_{\#}^k : G \to G$ that sends
each edge $E$ in $G$ to  the tight edge-path $f_{\#}^k(E)$ that is 
homotopic rel endpoints to $f^k(E)$.

\subsection{(Improved) Relative train tracks}

We now describe the properties of Improved Relative Train Track maps,
as constructed in \cite{BH2} and \cite{BFH}.

\bigskip
 \noindent{\bf{Splitting:}}
Suppose that $\sigma = \sigma_1\sigma_2$ is a decomposition of a path
 into nontrivial subpaths (we do not assume that
$\sigma_1$ and $\sigma_2$ are edge-paths, even if $\sigma$ is).  We
say that $\sigma = \sigma_1\sigma_2$ is a {\em $k$-splitting} if
\[      f_{\#}^k(\sigma) = f_{\#}^k(\sigma_1) f_{\#}^k(\sigma_2)
\]
is a decomposition into sub-paths (i.e. for {\em some} choice of tightening,
there is no folding between the $f$-images of $\sigma_1$ and
$\sigma_2$ for at least $k$ iterates). If $\sigma = \sigma_1\sigma_2$
is a $k$-splitting for all $k > 0$ then it is called a {\em splitting\footnote{In the next
  section, we introduce a stronger notion of {\em hard} splittings.}}
and we write $\sigma = \sigma_1 \cdot \sigma_2$.
If one of $\sigma_1$
or $\sigma_2$ is the empty path, the splitting is said to be {\em
  vacuous}.

A {\em turn} in $G$ is an unordered pair of half-edges originating at
a common vertex.  A turn is {\em non-degenerate} if it is defined by
distinct  half-edges, and is {\em degenerate} otherwise.  The
map $f : G \to G$ induces a self-map $Df$ on the set of oriented edges
of $G$ by sending an oriented edge to the first oriented edge in its
$f$-image.   $Df$ induces a map $Tf$ on the set of turns in $G$.

A turn is {\em illegal} with respect to $f : G \to G$ if its image
under some iterate of $Tf$ is degenerate; a turn is {\em legal} if it
is not illegal.

Associated to  $f$ is a {\em filtration} of $G$, 
\[      \emptyset = G_0 \subset G_1 \subset \cdots \subset G_\omega= G,   \]
consisting of $f$-invariant subgraphs of $G$.  We call the sets $H_r
:= \overline{ G_r \ssm G_{r-1}}$ {\em strata}. To each stratum $H_r$
is associated $M_r$, the {\em transition matrix for $H_r$}; the
$(i,j)^{\rm th}$ entry of $M_r$ is the number of times the $f$-image
of the $j^{\rm th}$ edge crosses the $i^{\rm th}$ edge in either
direction. By choosing a filtration carefully one may ensure that for
each $r$ the matrix $M_r$ is either the zero matrix or is 
irreducible. If
$M_r$ is the zero matrix, then we say that $H_r$ is a {\em zero
stratum}.  Otherwise, $M_r$ has an associated Perron-Frobenius
eigenvalue $\lambda_r \geq 1$, see \cite{Seneta}. If $\lambda_r > 1$
then we say that $H_r$ is an {\em exponential stratum}; if $\lambda_r
= 1$ then we say that $H_r$ is a {\em parabolic stratum}\footnote{Bestvina {\em et al.} use the terminology {\em
exponentially-growing} and {\em non-exponentially-growing} for our
exponential and parabolic. This difference in terminology explains the
names of the items in Theorem \ref{MainTrainTrack} below.}. The edges
in strata inherit these adjectives, e.g. ``exponential edge".  

A turn is defined to be  {\em in $H_r$ if both} half-edges 
lie in the stratum $H_r$.  A turn is a {\em
mixed turn in $(G_r,G_{r-1})$} if one edge is in $H_r$ and the other is in
$G_{r-1}$. A path with no illegal turns in $H_r$ is said to be {\em
$r$-legal}.  We may emphasize that certain turns are in $H_r$ by calling them 
{\em $r$-(il)legal turns}.

\begin{definition} \cite[Section 5, p.38]{BH2} \label{TrainTrackDef}
We say that $f : G \to G$ is a {\em relative train track map} if the
following conditions hold for every exponential stratum $H_r$: 
\begin{enumerate}
\item[\rm (RTT-i)] $Df$ maps the set of oriented edges in $H_r$ to
  itself; in particular all mixed turns in $(G_r,G_{r-1})$ are legal. 
\item[\rm (RTT-ii)] If $\alpha$ is a nontrivial path in $G_{r-1}$
  with endpoints in $H_r \cap G_{r-1}$, then $f_{\#}(\alpha)$ is a
  nontrivial path with endpoints in $H_r \cap G_{r-1}$. 
\item[\rm (RTT-iii)] For each legal path $\beta$ in $H_r$,
  $f(\beta)$ is a path that does not contain any illegal turns in
  $H_r$. 
\end{enumerate}
\end{definition}

The following lemma is ``the most important consequence of being a
relative train track map'' \cite[p.530]{BFH}; it follows immediately
from Definition \ref{TrainTrackDef}. 

\begin{lemma} \cite[Lemma 5.8, p.39]{BH2} \label{RTT(i)Lemma}
Suppose that $f : G \to G$ is a relative train track map, that $H_r$
is an exponential stratum and that $\sigma = a_1b_1a_2 \dots b_l$ is
the decomposition of an $r$-legal path $\sigma$ into subpaths $a_j$ 
in $H_r$ and $b_j$ in $G_{r-1}$. (Allow for the possibility
that $a_1$ or $b_l$ is trivial, but assume the other subpaths are
nontrivial.) Then $f_{\#}(\sigma) = f(a_1)f_{\#}(b_1)f(a_2)\dots
f_{\#}(b_l)$ and is $r$-legal.
\end{lemma}

\begin{definition} \label{BasicPaths}
Suppose that $f : G \to G$ is a topological representative, that the
parabolic stratum $H_i$ consists of a single edge
$E_i$ and that $f(E_i) = E_iu_i$ for some path $u_i$ in $G_{i-1}$.
We say that the paths of the form $E_i\gamma \bar{E_i}$, $E_i\gamma$
and $\gamma \bar{E_i}$, where $\gamma$ is in $G_{i-1}$, are {\em basic paths
of height $i$}.
\end{definition}

\begin{lemma} \cite[Lemma 4.1.4, p.555]{BFH}
\label{BasicPathSplitting}
Suppose that $f : G \to G$ and $E_i$ are as in Definition
\ref{BasicPaths}.  Suppose further that $\sigma$ is a
path or circuit in $G_i$ that intersects $H_i$ nontrivially and that the
endpoints of $\sigma$ are not contained in the interior of $E_i$.
Then $\sigma$ has a splitting each of whose pieces is either a  basic
path of height $i$ or is contained in $G_{i-1}$.
\end{lemma}

\begin{definition}
A {\em Nielsen path} is a nontrivial path $\sigma$ such that
$f^k_{\#}(\sigma) = \sigma$ for some $k \geq 1$.
\end{definition}
Nielsen paths are called {\em periodic Nielsen paths} in \cite{BFH}, but
Theorem \ref{MainTrainTrack} below allows us to choose an $f$ so  that
any periodic Nielsen path has period $1$ (which is to say that
$f_{\#}(\sigma) = \sigma$), and we shall assume that $f$ satisfies the
properties outlined in Theorem \ref{MainTrainTrack}.  Thus we can
assume that $k=1$ in the above definition.  A Nielsen path is called
{\em indivisible} if it cannot be split as a concatenation of two
non-trivial Nielsen paths.

\begin{definition} [cf. 5.1.3, p.531 \cite{BFH}] \label{Exceptional}
Suppose that $H_i$ is a single edge $E_i$ and that $f(E_i) =
E_i\tau^l$ for some closed Nielsen path $\tau$ in
$G_{i-1}$ and some $l > 0$.  The {\em exceptional paths of height $i$}
are those paths of the form $E_i\tau^k\bar{E}_j$ or
$E_i\bar{\tau}^k\bar{E}_j$ where $k \geq 0$, $j \leq i$, $H_j$ is a
single edge $E_j$ and $f(E_j) = E_j\tau^m$ for some $m > 0$.
\end{definition}

\begin{remark} In \cite{BFH} the authors require that the path $\tau$ is an
  indivisible Nielsen path.  However, exceptional paths are defined so
  that condition ne-(iii) of Theorem \ref{MainTrainTrack} holds, and
  an analysis of the proof of this theorem in \cite{BFH} shows that
  the restriction to indivisible Nielsen paths in exceptional paths is
  invalid.
\end{remark}

In Definition \ref{Exceptional}, the paths do not have a preferred
orientation.  Thus it is
important to note that the paths of the form $E_j \tau^k\bar{E_i}$ and
$E_j \bar{\tau}^k\bar{E_i}$ with $E_i,E_j$ and $\tau$ as above are
also exceptional paths of height $i$.

\subsection{The Theorem of Bestvina, Feighn and Handel}

A matrix is {\em aperiodic} if it has a power in which every entry is
positive.  The map $f$ is {\em eg-aperiodic} if every
exponential stratum has an aperiodic transition matrix.

Theorem 5.1.5 in \cite{BFH} is the main structural theorem for
improved relative train track maps.  We shall use it continually in
what follows, often without explicit mention.  We therefore record
those parts of it which we need. A map $f$ which satisfies the
statements of Theorem \ref{MainTrainTrack} is called an {\em improved
relative train track map}.

\begin{theorem} (cf. Theorem 5.1.5, p.562, \cite{BFH}) \label{Theorem5.1.5}
  \label{MainTrainTrack} 
For every outer automorphism ${\mathcal O} \in {\rm Out}(F)$ there is an
eg-aperiodic relative train track map $f : G \to G$ with filtration
$\emptyset = G_0 \subset G_1 \subset \dots \subset G_\omega = G$ such that
$f$ represents an iterate of ${\mathcal O}$, and $f$ has the following
properties.
\begin{itemize}
\item Every periodic Nielsen path has period one.
\item For every vertex $v\in G$, $f(v)$ is a fixed point.  If $v$ is
  an endpoint of an edge in a parabolic stratum then $v$ is a fixed
  point.  If $v$ is the endpoint of an edge in an
  exponential stratum $H_i$ and if $v$ is also contained in
  a noncontractible component of $G_{i-1}$, then $v$ is a fixed
  point. 
\item $H_i$ is a zero stratum if and only if it is the union of the
  contractible components of $G_i$. 
\item If $H_i$ is a {\em zero stratum}, then
\begin{enumerate}
\item[\rm z-(i)] $H_{i+1}$ is an exponential stratum.
\item[\rm z-(ii)] $f|H_i$ is an immersion.
\end{enumerate}
\item If $H_i$ is a {\em parabolic stratum}, then
\begin{enumerate}
\item[\rm ne-(i)] $H_i$ is a single edge $E_i$.
\item[\rm ne-(ii)] $f(E_i)$ splits as $E_i \cdot u_i$ for some closed
  path $u_i$ in $G_{i-1}$ whose basepoint is fixed by $f$.
\item[\rm ne-(iii)] If $\sigma$ is a basic path of height
  $i$ that does not split as a concatenation of two basic paths of
  height $i$ or as a concatenation of a basic path of height $i$ with
  a path contained in $G_{i-1}$, then either: {\rm (i)} for some $k$,
  the path $f^k_{\#}(\sigma)$ splits  into pieces, one of which equals
  $E_i$ or $\bar{E_i}$; or {\rm (ii)} $u_i$ is a Nielsen path and, for
  some $k$, the path $f^k_{\#}(\sigma)$ is an exceptional path of
  height $i$. 
\end{enumerate}
\item If $H_i$ is an {\em exponential stratum} then
\begin{enumerate}
\item[\rm eg-(i)] There is at most one indivisible Nielsen path
  $\rho_i$ in $G_i$ that intersects $H_i$ nontrivially.  The
  initial edges of $\rho_i$ and $\bar{\rho_i}$ are distinct (possibly
  partial) edges in $H_i$.
\end{enumerate}
\end{itemize}
\end{theorem}

Suppose that $f : G \to G$ is an improved relative train track map
representing some iterate $\phi^k$ of $\phi \in \On$, and that $\rho$
is a Nielsen path in $G_r$ that intersects $H_r$ nontrivially,
and suppose that $\rho$ is not an edge-path.  Then subdividing the
edges containing the endpoints of $\rho$ at the endpoints, gives a new
graph $G'$, and the map $f' : G' \to G'$ induced by $f$ is an improved
relative train track map representing $\phi^k$.  To ease notation, it
is convenient to assume that this subdivision has been performed.
Under this assumption, all Nielsen paths will be edge-paths, and all
of the paths which we consider in the remainder of this paper will
also be edge-paths.

\begin{convention}
{\rm Since all Nielsen paths in the remainder of this paper will be edge paths, we will use the phrase {\em `indivisible Nielsen path'} to mean a Nielsen edge-path which cannot be decomposed nontrivially as a concatenation of two non-trivial Nielsen {\em edge}-paths.  In particular, a single edge fixed pointwise by $f$ will be considered to be an indivisible Nielsen path.}
\end{convention} 

\medskip

{\em For the remainder of this article, we will concentrate on an improved relative train track map $f : G \to G$ and repeatedly pass to iterates $f_{\#}^k$ in order to better control its cancellation properties.}

\medskip

Recall the following from \cite[Section 4.2, pp.558-559]{BFH}.

\begin{definition}
If $f:G \to G$ is a relative train track map and $H_r$ is an
exponential stratum, then define $P_r$ to be the set of paths $\rho$
in $G_r$ that are such that: 
\begin{enumerate}
\item[(i)] For each $k \geq 1$ the path $f^k_{\#}(\rho)$ contains
  exactly one illegal turn in $H_r$. 
\item[(ii)] For each $k \geq 1$ the initial and terminal (possibly
  partial) edges of $f^k_{\#}(\rho)$ are contained in $H_r$.
\item[(iii)] The number of $H_r$-edges in $f^k_{\#}(\rho)$ is bounded
  independently of $k$.
\end{enumerate}
\end{definition}

\begin{lemma} \cite[Lemma 4.2.5, p.558]{BFH} \label{Pr}
$P_r$ is a finite $f_{\#}$-invariant set.
\end{lemma}

\begin{lemma} \cite[Lemma 4.2.6, p.559]{BFH} \label{ExpSplitting}
Suppose that $f : G \to G$ is a relative train track map, that $H_r$
is an exponential stratum, that $\sigma$ is a
path or circuit in $G_r$ and that, for each $k\geq 0$, the path
$f^k_{\#}(\sigma)$ has the same finite number of illegal turns in
$H_r$.  Then $\sigma$ can be split into subpaths that are either
$r$-legal or elements of $P_r$.
\end{lemma}

\begin{definition} \label{Weightr}
If $\rho$ is a path and $r$ is the least integer such that $\rho$ is in 
$G_r$ then we say that {\em $\rho$ has weight $r$}.
\end{definition}

If $\rho$ has weight $r$ and $H_r$ is exponential, we will say that
$\rho$ is an {\em exponential path}.  We define {\em
parabolic paths} similarly.

\begin{lemma} \label{PreNielsen}
Suppose that $\sigma$ is an edge-path and that, for some $k \geq 1$,
$f_{\#}^k(\sigma)$ is a Nielsen path.  Then $f_{\#}(\sigma)$ is a
Nielsen path. 
\end{lemma}

\begin{proof} Suppose that the endpoints of $\sigma$ are $u_1$ and
$v_1$ and that the endpoints of $f_{\#}^k(\sigma)$ are $u_2$ and
$v_2$.  For each vertex $v \in G$, $f(v)$ is fixed by $f$, so
$f(u_1) = u_2$ and $f(v_1) = v_2$.  If $f_{\#}(\sigma) \neq
f_{\#}^k(\sigma)$ then we have two edge-paths with the same
endpoints which eventually get mapped to the same path.  Thus there is
some nontrivial circuit which is killed by $f$, contradicting the fact
that $f$ is a homotopy equivalence.  Therefore $f_{\#}(\sigma) =
f_{\#}^k(\sigma)$ and so is a Nielsen path.
\end{proof}

\medskip

Always, $L$ will denote the maximum of $|f(E)|$, for $E$
an edge in $G$.

Later, we will pass to further iterates of $f$ in order to find a
particularly nice form.

An analysis of the results in this section allows us to see that there
are three kinds of irreducible Nielsen paths.  The first are those
which are single edges; the second are certain exceptional paths; and the
third lie in the set $P_r$.  We will use this trichotomy frequently
without mention.  The first two cases are where the path is
parabolic-weight, the third where it is exponential-weight.  It is not
possible for Nielsen path to have weight $r$ where $H_r$ is a zero
stratum. 

\begin{observation}
Let $\rho$ be an indivisible Nielsen path of weight $r$.  Then the
first and last edges in $\rho$ are contained in $H_r$.
\end{observation}

Because periodic Nielsen paths have period $1$, the set of Nielsen
paths does not change when $f$ is replaced by a further iterate of
itself.  We will use this fact often.

\begin{lemma} \label{LinearEdges}
Suppose $E$ is an edge such that $|f_{\#}^j(E)|$ grows linearly with
$j$.  Then $f(E) = E \cdot \tau^k$, where $\tau$ is a Nielsen path
that is not a proper power.  The edge-path $\tau$ decomposes into
indivisible Nielsen paths each of which is itself an edge-path. 
\end{lemma}
\begin{proof}
The fact that $f(E) = E \cdot \tau^k$, where $\tau$ is a Nielsen path
follows from conditions ne-(ii) and ne-(iii) of Theorem
\ref{MainTrainTrack}.  \footnote{If Theorem \ref{MainTrainTrack},
ne-(iii) held with the Nielsen path $\tau$ in the definition of
exceptional paths being indivisible, we could also insist that
$\tau$ be indivisible here.}

The final sentence follows because $\tau$ is an edge-path.  Any
Nielsen path admits a splitting into indivisible Nielsen paths.  If
there were an indivisible Nielsen path in this decomposition which was
not an edge-path then it
would have to be of exponential weight, and there is at most one
indivisible Nielsen path of each exponential weight.  Therefore, the
final Nielsen path of this weight would end at a half-edge, and the
remainder of this edge could not be contained in an indivisible
Nielsen path, contradicting the decomposition of $\tau$ into
indivisible Nielsen paths.
\end{proof}

\begin{lemma} \label{NoTau}
Let $\tau$ be a Nielsen path and $\tau_0$ a proper initial (or
terminal) sub edge-path of $\tau$.  No image $f_{\#}^k(\tau_0)$ contains
$\tau$ as a sub-path.
\end{lemma}
\begin{proof} 
It is sufficient to prove the lemma for indivisible Nielsen paths, as
the result for arbitrary Nielsen paths then follows immediately.

If $\tau$ is an indivisible Nielsen path and $\tau_0$ is a proper
non-trivial subpath of $\tau$ then $\tau$ cannot be a single edge.
Therefore, either $\tau$ is either an indivisible Nielsen path of
exponential weight, or an exceptional path.

In case $\tau$ is an indivisible Nielsen path of exponential weight,
suppose the weight is $r$.  Then, by Lemma \ref{ExpSplitting} $\tau$
contains a single illegal turn in $H_r$.  If $\tau_0$ does not contain
this illegal turn then $\tau_0$ is $r$-legal, and so no iterate of
$\tau_0$ contains an illegal turn in $H_r$, and therefore no iterate
of $\tau_0$ can contain $\tau$ as a subpath.  If $\tau_0$ does contain
the $r$-illegal turn in $\tau$ then, being a proper subpath of
$\tau$, the path on one side of the illegal turn in $\tau_0$ and its
iterates is strictly smaller than the corresponding path in $\tau$,
and again $\tau$ cannot be contained as a subpath of any iterate of
$\tau_0$. 

If $\tau$ is an exceptional path, then $\tau = E_i \rho^k
\bar{E_j}$ where $\rho$ is a Nielsen path and $E_i$ and $E_j$ are
of weight greater than $\rho$.  Therefore, any sub edge-path $\tau_0$
of $\tau$ contains at most one edge of weight greater than $\rho$, and
the same is true for any iterate of $\tau_0$, and once again no
iterate of $\tau_0$ contains $\tau$ as a sub-path. 
\end{proof}

\section{Hard splittings} \label{s:hard}

In this section we introduce a new concept for improved relative train tracks:
{\em hard splittings}.  This plays an important role in the subsequent sections
of this paper, and also in \cite{bg3}.

Recall that a decomposition of a path $\sigma = \sigma_1 \sigma_2$ is
a $k$-splitting if $f_{\#}^k(\sigma) = f_{\#}^k(\sigma_1)
f_{\#}^k(\sigma_2)$; which means that, for {\em some} choice of
tightening, the images of $\sigma_1$ and $\sigma_2$ do not interact
with each other.  This leads to the concept of {\em splittings}.  We
need a more restrictive notion, where the decomposition is preserved
for {\em every} choice of tightening. For this purpose, we make the
following

\begin{definition} [Hard splittings]
We say that a $k$-splitting $\rho = \rho_1 \rho_2$ is a {\em
  hard} $k$-splitting if for {\em any} choice of tightening of $f^k(\rho) =
f^k(\rho_1) f^k(\rho_2)$ there is no cancellation between the image of
$\rho_1$ and the image of $\rho_2$.

A decomposition which is a hard $k$-splitting for all $k \ge 1$ is
called a {\em hard} splitting.  If $\rho_1 \cdot \rho_2$ is a hard
splitting, we write $\rho_1 \odot \rho_2$.

An edge path is {\em hard-indivisible} (or {\em h-indivisible}) if it admits no non-vacuous hard splitting into edge paths.
\end{definition}

\begin{example}
Suppose that $G$ is the graph with a single vertex and edges $E_1,
E_2$ and $E_3$.  Suppose that $f(E_1) = E_1$, $f(E_2) = E_2 E_1$ and
$f(E_3) = E_3 \bar{E_1} \bar{E_2}$.  Then $f$ is an improved
relative train track.  Then $E_3E_2 \cdot \bar{E_1}$ is a
$1$-splitting, since
\[      f(E_3E_2\bar{E_1}) = E_3 \bar{E_1}\bar{E_2} E_2
E_1 \bar{E_1},      \]
which tightens to $E_3 \bar{E_1} = f_{\#}(E_3E_2)
f_{\#}(\bar{E_1})$. In fact this is a splitting. However, there is
a choice of tightening which first cancels the final $E_1
\bar{E_1}$ and then the subpath $\bar{E_2} E_2$.  Therefore
the splitting $E_3 E_2 \cdot \bar{E_1}$ is not a hard
$1$-splitting.
\end{example}

The following lemma shows the main utility of hard splittings, and the
example above shows that it is not true in general for splittings. 

\begin{lemma} \label{HardSplitLemma}
Suppose that $\sigma_1 \odot \sigma_2$ is a hard splitting, and that
$\rho$ is an initial subpath of $\sigma_2$.  Then $\sigma_1 \odot
\rho$ is a hard splitting. 
\end{lemma}
\begin{proof}
If there were any cancellation between
images of $\sigma_1$ and $\rho$ then there would be a possible
tightening between the images of $\sigma_1$ and $\sigma_2$. 
\end{proof}

The following two lemmas will also be crucial for our applications of
hard splittings in \cite{bg3}.

\begin{lemma}
Every edge path admits a unique maximal hard splitting into edge paths.
\end{lemma}
\begin{proof} This follows by an obvious induction on length from the
observation that
if $\rho = \rho_1 \rho_2 \rho_3$,
where the $\rho_i$ are edge paths,  and if
$\rho = \rho_1 \odot \rho_2 \rho_3$ and $\rho = \rho_1 \rho_2 \odot \rho_3$
then $\rho = \rho_1 \odot \rho_2 \odot \rho_3$. 
\end{proof}

\begin{lemma}
If $\rho = \rho_1 \odot \rho_2$ and $\sigma_1$ and $\sigma_2$
are, respectively, terminal and initial subpaths of $f^k_{\#}(\rho_1)$
and $f^k_{\#}(\rho_2)$ for some $k \ge 0$ then
$\sigma_1 \sigma_2 = \sigma_1 \odot \sigma_2$.
\end{lemma}
\begin{proof}
For all $i \ge 0$, the untightened path $f^i(\sigma_1)$
is a terminal subpath of the untightened path $f^{k+i}(\rho_1)$,
while $f^i(\sigma_2)$ is an initial subpath of $f^{k+i}(\rho_2)$.

The hardness of the splitting $\rho = \rho_1 \odot \rho_2$ ensures
that no matter how one tightens $f^{k+i}(\rho_1)f^{k+i}(\rho_2)$
there will be no cancellation between $f^{k+i}(\rho_1)$ and $f^{k+i}(\rho_2)$.
In particular, one is free to 
tighten $f^i(\sigma_1)f^i(\sigma_2)$ first, and 
there can be no cancellation between them.
(It may happen that when one goes to tighten $f^{k+i}(\rho_1)$ completely,
the whole of $f^i(\sigma_1)$ is cancelled, but this does not affect the assertion
of the lemma.)
\end{proof}

The purpose of the remainder of this section is to sharpen results
from the previous section to cover hard splittings. \footnote{Bestvina 
{\em et al.} make no explicit mention of the
  distinction between splittings and hard splittings, however
  condition (3) of Proposition 5.4.3 on p.581 (see Lemma
  \ref{SubpathSplitting} below) indicates that they are aware of the
  distinction and that the term `splitting' has the same
  meaning for them as it does here.}

\begin{lemma} [cf. Lemma 4.1.1, p.554 \cite{BFH}] \label{HardSplitProperties}
If $\sigma = \sigma_1 \odot \sigma_2$ is a hard splitting, and
  $\sigma_1 = \sigma_1' \odot \sigma_2'$ is a hard splitting then
  $\sigma = \sigma_1' \odot \sigma_2' \odot \sigma_2$ is a hard
  splitting.  The analogous result with the roles of $\sigma_1$ and
  $\sigma_2$ reversed also holds.
\end{lemma}

\begin{remark} The possible existence of an edge-path $\sigma_2$ so that
$f_{\#}(\sigma_2)$ is a single vertex means that
 $\sigma_1 \sigma_2 = \sigma_1 \odot \sigma_2$ and $\sigma_2 \sigma_3 = \sigma_2 \odot \sigma_3$ need {\em not} imply
that  $\sigma_1 \sigma_2 \sigma_3 = \sigma_1 \odot \sigma_2 \odot \sigma_3$.

Indeed if $\sigma_2$ is an edge-path so that 
 $f_{\#}(\sigma_2)$ is a vertex then
$f_{\#}(\sigma_1)$ and $f_{\#}(\sigma_3)$ come together in a tightening of
$f(\sigma_1\sigma_2\sigma_3)$, possibly cancelling.

In contrast, if  $f_{\#}(\sigma_2)$ (and hence each $f_{\#}^k(\sigma_2)$)
contains an edge, then  the hardness of the two splittings ensures that
in any tightening $f_{\#}(\sigma_1\sigma_2\sigma_3)=f_{\#}(\sigma_1)f_{\#}(\sigma_2)
f_{\#}(\sigma_3) $, that is $\sigma_1 \sigma_2 \sigma_3 = \sigma_1 \odot \sigma_2 \odot \sigma_3$.
\end{remark}

The following strengthening of Theorem \ref{MainTrainTrack} ne-(ii) is
a restatement of (a weak form of) \cite[Proposition 5.4.3.(3),
  p.581]{BFH}.

\begin{lemma} \label{SubpathSplitting}
Suppose $f$ is an improved relative train track map and $E$ is a parabolic edge with $f(E) = E u$.  For any initial subpath $w$ of $u$, $E \cdot w$ is
a splitting.
\end{lemma}

\begin{corollary}
Suppose $f$ is an improved relative train track map, $E$ is a
parabolic edge and $f(E) = E u$.  Then $f(E) = E \odot u$.
\end{corollary}

\begin{lemma} \label{BasicHSplit}
Suppose $H_i$ is a parabolic stratum and $\sigma$ is a
path in $G_i$ that intersects $H_i$ nontrivially, and that the endpoints of
$\sigma$ are not contained in the interior of $E_i$.  Then $\sigma$
admits a hard splitting, each of whose pieces is either a basic path
of height $i$ or is contained in $G_{i-1}$. 
\end{lemma}

\begin{lemma} \label{ne(iii)Hard}
If $\sigma$ is a basic path of height $i$ that does not
admit a hard splitting as a concatenation of two basic paths of height 
$i$ or as a concatenation of a basic path of height $i$ with a path of 
weight less than $i$, then either; (i) for some $k$, the path
$f_{\#}^k(\sigma)$ admits a hard splitting into pieces, one of which
is $E_i$ or $\bar{E_i}$; or (ii) $f(E_i) = E_i \odot u_i$, where
$u_i$ is a Nielsen path and, for some $k$, the path $f_{\#}^k(\sigma)$
is an exceptional path of height $i$.
\end{lemma}

\begin{proof}
Follows from the proof of \cite[Lemma 5.5.1, pp.585--590]{BFH}.
\end{proof}

\begin{lemma} [cf. Lemma \ref{ExpSplitting} above] \label{ExpHardSplitting}
Suppose that $f : G \to G$ is a relative train track map, that $H_r$
is an exponentially-growing stratum, that $\sigma$ is a
path or circuit in  $G_r$, and that each $f_{\#}^k(\sigma)$ has the same finite
number of illegal turns in $H_r$.  Then $\sigma$ can be decomposed as
$\sigma = \rho_1 \odot \ldots \odot \rho_k$, where
each $\rho_i$ is either (i) an element of $P_r$;
(ii) an $r$-legal path which starts and ends with edges in $H_r$; or
(iii) of weight at most $r-1$.
\end{lemma}
\begin{proof}
Consider the splitting of $\sigma$ given by Lemma \ref{ExpSplitting}.
The pieces of this splitting are either (i) elements of $P_r$, or (ii)
$r$-legal paths.  By Definition \ref{TrainTrackDef} RTT-(i), any
$r$-legal path admits a hard splitting into $r$-legal paths which
start and end with edges in $H_r$, and paths of weight at most $r-1$.
The turn at the end of a Nielsen path in the splitting of $\sigma$ is
either a mixed turn (with the edge from $H_r$ coming from the Nielsen
path and the other edge being of weight at most $r-1$) or a legal turn
in $H_r$.  In either case, $\sigma$ admits a hard splitting at the
vertex of this turn.
\end{proof}

The next result follows from a consideration of the form of indivisible Nielsen paths, noting Definition \ref{TrainTrackDef} and Lemma \ref{ExpHardSplitting}.

\begin{lemma} \label{NielsenHardSplit}
Any Nielsen path admits a hard splitting into indivisible Nielsen paths.
\end{lemma}

\begin{remark} \label{HardSplitRemark}
If $\rho = \rho_1 \odot \rho_2$ is a hard splitting for the map $f$
then it is a hard splitting for $f_{\#}^k$ for any $k \ge 1$.
\end{remark}

We record a piece of terminology for later use.

\begin{definition} \label{Displayed}
A sub edge-path $\rho$ of a path $\chi$ is {\em displayed} if there is a hard splitting of $\chi$ immediately on either side of $\rho$.
\end{definition}

 \section{A small reduction}
\label{ColourCancellation}

In this section we clarify a couple of issues about monochromatic paths,
and state Theorem \ref{MainProp}, which immediately implies the Beaded
Decomposition Theorem.

Our strategy for proving the Beaded Decomposition Theorem is as follows:
given an automorphism $\phi \in$\Aut, we start
with an improved relative train track representative $f :G \to G$ for some iterate $\phi^k$ of $\phi$, as obtained from the conclusion of Theorem \ref{MainTrainTrack}.  We analyse
the evolution of monochromatic paths, and eventually pass to an iterate of $f$
in which we can prove the Beaded Decomposition Theorem.  However, it is crucial
to note that monochromatic paths for $f$ are not necessarily monochromatic paths
for $f^k_{\#}$ when $k > 1$.  See Section \ref{IterateSection} for further disussion about some of these issues.  

These concerns lead to the following definition, where we are
concentrating on a fixed IRTT $f : G \to G$, and so omit mention of $f$ from our notation.

\begin{definition} \label{MonoChiDef}
For a positive integer $r$, we define {\em $r$-monochromatic paths}  by  recursion:
edges in $G$ are $r$-monochromatic and if $\rho$ is an
$r$-monochromatic path then every sub edge-path of
$f_\#^{r}(\rho)$ is $r$-monochromatic.
\end{definition}

Note that if $r'$ is a multiple of $r$ then every $r'$-monochromatic path is  $r$-monochromatic
but not {\em vice versa}. Thus if we replace $f$ by an iterate then, for fixed $n$, 
the set of 
$n$-monochromatic paths may get  smaller. The content of the Beaded Decomposition Theorem is that
one need only pass to a bounded iterate in order to ensure that all monochromatic
paths admit a beaded decomposition.  In particular, the Beaded Decomposition
Theorem is an immediate consequence of the following theorem.

\begin{theorem} [Monochromatic paths are beaded] \label{MainProp}
Let $f:G\to G$ be an improved relative train track map.
There exist constants $r$ and $J$, depending only on $f$,
so that every $r$-monochromatic path in $G$ is $(J,f)$-beaded.
\end{theorem}

\section{Nibbled futures}\label{s:nib}

Monochromatic paths arise as  {\em nibbled futures} in the sense defined
below. Thus   in order to prove 
Theorem \ref{MainProp} we must understand how {\em nibbled futures} evolve.
The results in this section reduce this challenge to the task of understanding
the nibbled futures of GEPs.

\begin{definition}[\NF s]
Let $\rho$ be a (tight) edge path.  The {\em $0$-step \nf \ of $\rho$} is $\rho$.

For $k \ge 1$, a {\em $k$-step \nf \ of $\rho$} is a sub edge-path of $f_{\#}(\sigma)$, where $\sigma$ is a $(k-1)$-step \nf \ of $\rho$.  A {\em \nf \ of $\rho$} is a $k$-step \nf \ for some $k \ge 0$.

For $k \ge 0$, the {\em $k$-step entire future of $\rho$} is $f_{\#}^k(\rho)$.
\end{definition}

\begin{remark} The $1$-monochromatic paths are precisely
 the \nf s of  single edges.
\end{remark}

\begin{theorem} [First Decomposition Theorem] \label{DecompTheorem}
For any $n \ge 1$ there exists an integer $V = V(n,f)$ such that if
$\rho$ is an edge path of length at most $n$ then any \nf \ of $\rho$
admits a hard splitting into edge paths, each of which is either the
\nf \ of a \gep \ or else has length at most $V$. 
\end{theorem}

The remainder of this section is dedicated to proving Theorem \ref{DecompTheorem}.
We begin by examining the entire future of a path of fixed length (Lemma 
\ref{SplittingLemma}) and then refine the argument to deal with nibbling. In
the proof of the first of these lemmas we require the following observation.
 
\begin{remark} \label{NoOfTurns}
Suppose that $\rho$ is a tight path of weight $r$.  The immediate past of an $r$-illegal turn in $f_{\#}(\rho)$ (under any choice of tightening) is an $r$-illegal turn, and two $r$-illegal turns cannot have the same $r$-illegal turn as their past. In particular, the number of $r$-illegal turns in $f_{\#}^l(\rho)$ is a non-increasing function of $l$, bounded below by $0$.
\end{remark}

\begin{lemma} \label{SplittingLemma}
There is a function $D :  \N \to \N$, depending only on $f$, such that, for any $r \in \{ 1 , \ldots , \omega \}$, if $\rho$ is a path of weight $r$, and $|\rho| \le n$,
then for any $i \ge D(n)$ the edge path $f_{\#}^{i}(\rho)$ admits a hard splitting into edge paths, each of
which is either 
\begin{enumerate}
\item   a single edge of weight $r$;
\item an indivisible Nielsen path of weight $r$;
\item a \gep \ of weight $r$; or
\item a path of weight at most $r-1$.
\end{enumerate}
\end{lemma}
\begin{proof}  If $H_r$ is a zero stratum, then $f_{\#}(\rho)$ has
  weight at most $r-1$, and $D(n) = 1$ will suffice for any $n$. 

If $H_r$ is a parabolic stratum, then $\rho$ admits a hard splitting
into pieces which are either basic of height $r$ or of weight at most
$r-1$ (Lemma \ref{BasicHSplit}).  Thus it is sufficient to consider the case where $\rho$ is a
basic path of weight $r$ and $|\rho| \le n$.  By at most $2$ applications
of Lemma \ref{ne(iii)Hard}, we see that there exists a $k$ such
that $f_{\#}^k(\rho)$ admits a hard splitting into pieces which are
either (i) single edges of weight $r$, (ii) exceptional paths of
height $r$, or (iii) of weight at most $r-1$.  By taking the maximum
of such $k$ over all basic paths of height $r$ which are of length at
most $n$, we find an integer $k_0$ so that we have the desired hard splitting 
of $f_{\#}^{k_0}(\rho)$ for all basic paths of height $r$ of length
at most $n$.  Any of the exceptional paths in these splittings which
are not \gep s have bounded length and are either indivisible Nielsen
paths or are decreasing in length.  A crude bound on the length of the exceptional paths which are not \gep s is $L^{k_0}n$ where $L$ is the maximum length of $f(E)$ over all edges $E
\in G$.  Thus, those exceptional paths which are decreasing in length will
become \gep s within less than $L^{k_0}n$ iterations.  Therefore, 
replacing $k_0$ by $k_0 + L^{k_0}n$, we may assume all exceptional
paths in the hard splitting are \gep s. 

Finally, suppose that $H_r$ is an exponential stratum.  As noted in Remark \ref{NoOfTurns}, the number of $r$-illegal turns in $f_{\#}^l(\rho)$ is a non-increasing function
of $l$ bounded below by $0$.  Therefore, there is some $j$ so that the
number of $r$-illegal turns in $f_{\#}^{j'}(\rho)$ is the same for
all $j' \ge j$.  By Lemma \ref{ExpHardSplitting},
$f_{\#}^j(\rho)$ admits a hard splitting into pieces which are either
(i) elements of $P_r$, (ii) single edges in $H_r$, or (iii) paths of weight 
at most $r-1$.  To finish the proof of the lemma it remains to note that 
if $\sigma \in P_r$ then $f_{\#}(\sigma)$ is a Nielsen path by Lemma \ref{PreNielsen}. 

Therefore, the required constant for $H_r$ may be taken to be the maximum of $j+1$ over all the
paths of weight $r$ of length at most $n$. 

To find $D(n)$ we need merely take the maximum of the constants found
above over all of the strata $H_r$ of $G$. 
\end{proof}

In the extension of the above proof to cover nibbled futures, we shall need
the following straightforward adaptation of Lemma \ref{NoTau}.

\begin{lemma} \label{NibNoTau}
Let $\tau$ be a Nielsen path and $\tau_0$ a proper initial (or terminal) sub-path of $\tau$.  No nibbled future of $\tau_0$ contains $\tau$ as a sub-path.
\end{lemma}

\begin{proposition} \label{NibbleSplittingLemma}
There exists a function $D' :  \N \to \N$, depending only on $f$, so that for any $r \in \{ 1 , \ldots , \omega \}$, if $\rho$ is a path of weight $r$ and $|\rho| \le n$, then for any $i \ge D'(n)$ any $i$-step \nf \ of $\rho$ admits a hard splitting into edge paths, each of which is either
\begin{enumerate}
\item   a single edge of weight $r$;
\item   a \nf \ of a weight $r$ indivisible Nielsen path;
\item   a \nf \ of a weight $r$ \gep ; or
\item   a path of weight at most $r-1$.
\end{enumerate}
\end{proposition}

\begin{remark} \label{stable}
Each of the conditions (1) -- (4) stated above is stable in the following sense: once an edge in a $k$-step \nf \ lies in a path satisfying one of these conditions, then any future of this edge in any further \nf \ will also lie in such a path (possibly the future will go from case (1) to case (4), but otherwise which case it falls into is also stable).  Thus we can split the proof of Proposition \ref{NibbleSplittingLemma} into a number of cases, deal with the cases separately by finding some constant which suffices, and finally take a maximum to find $D'(n)$.  An entirely similar remark applies to a number of subsequent proofs, in particular Theorem \ref{ColourCancelMain}.
\end{remark}

\begin{proof}[Proof (Proposition \ref{NibbleSplittingLemma}).]
Let $\rho_0 = \rho$ and for $j > 0$ let $\rho_j$ be a sub edge-path of $f_{\#}(\rho_{j-1})$.

If $H_r$ is a zero stratum, then $f_{\#}(\rho)$ has weight at most $r-1$ and it suffices to take $D'(n) = 1$.

Suppose that $H_r$ is an exponential stratum.  By Lemma \ref{SplittingLemma}, the $D(n)$-step entire future of $\rho$ admits a hard splitting of the desired form.  We consider how \nib \ can affect this splitting.  As we move forwards through the \nf \ of $\rho$, cancellation of $H_r$-edges can occur only at $r$-illegal turns and at the ends, where the \nib \ occurs. 

Remark \ref{NoOfTurns} implies that we can trace the $r$-illegal turns forwards through the successive \nf s of $\rho$ (whilst the $r$-illegal continues to exist).  We compare the $r$-illegal turns in $\rho_k$ to those in $f_{\#}^k(\rho)$, the entire future of $\rho$.  We say that the \nib \ {\em first cancels an $r$-illegal turn at time $k$} if the collection of $r$-illegal turns in $\rho_{k-1}$ is the same as the collection in $f_{\#}^{k-1}(\rho)$, but the collection in $\rho_k$ is {\em not} the same as that of $f_{\#}^k(\rho)$.  The first observation we make is that if, at time $k$, the \nib \ has not yet cancelled any $r$-illegal turn then the sequence of $H_r$-edges $\rho_k$ is a subsequence of the $H_r$-edges in $f_{\#}^k(\rho)$.  Therefore, any splitting of the desired type for $f_{\#}^k(\rho)$ is inherited by $\rho_k$.

Since there is a splitting of the $D(n)$-step entire future of $\rho$ of the desired form, either there is a splitting of $\rho_{D(n)}$, or else $\rho_{D(n)}$ has fewer $r$-illegal turns than $f_{\#}^{D(n)}(\rho)$, and hence than $\rho$.  However, $|\rho_{D(n)}| \le n.L^{D(n)}$.  We apply the above argument to $\rho_{D(n)}$, going forwards a further $D(nL^{D(n)})$ steps into the future.  Since the number of illegal turns in $H_r$ in $\rho$ was at most $n-1$, we will eventually find a splitting of the required form within an amount of time bounded by a function of $n$ (this function depends only on $f$, as required).  Denoting this function by $D_0$, we have that any $D_0(n)$-step \nf \ of any path of exponential weight whose length is at most $n$ admits a hard splitting of the desired form.

Now suppose that $H_r$ is a parabolic stratum.  By Lemma \ref{BasicHSplit}, $\rho$ admits a hard splitting into basic edge paths.  Therefore we may assume (by reversing the orientation of $\rho$ if necessary) that $\rho = E_r \sigma$ or $\rho = E_r \sigma \overline{E_r}$ where $E_r$ is the unique edge in $H_r$ and $\sigma$ is in $G_{r-1}$.  For the \nf \ of $\rho$ to have weight $r$, the \nib \ must occur only on one side (since the only edges of weight $r$ in any future of $\rho$ occur on the ends).  We assume that all \nib \ occurs from the right.  Once again, the $D(n)$-step entire future of $\rho$ admits a hard splitting of the desired form.  If $\rho = E_r \sigma \overline{E_r}$ then the $D(n)$-step \nf \ of $\rho$ either admits a hard splitting of the required form, or is of the form $E_r \sigma_1$, where $\sigma_1$ is in $G_{r-1}$.  Hence we may assume that $\rho = E_r \sigma$.  Suppose that $f(E_r) = E_r u_r$, and that $u_r$ has weight $s < r$.

We first consider the possibility that $f_{\#}(\sigma)$ has weight $q > s$ (but less than $r$ by hypothesis).  There are two cases to consider here.  The first is that $H_q$ is an exponential stratum.  The future of $E_r$ cannot cancel any edges of weight $q$ or higher in the future of $\sigma$, so the edges of weight $q$ in the \nf \ of $\rho$ are exactly the same as the edges of weight $q$ in the corresponding \nf \ of $\sigma$ (recall we are assuming that \nib \ only occurs from the right).  This $D_0(|\sigma|)$-step \nf \ of $\sigma$ admits a hard splitting into edge paths which are either\footnote{\gep s have parabolic weight} single edges of weight $q$, the \nf \ of an indivisible Nielsen path of weight $q$, or of weight at most $q-1$.  Let $\sigma_2$ be the path from the right endpoint of $E_r$ up to but not including the first edge of weight $q$.  Then, since mixed turns are legal, the $D_0(n)$-step \nf \ of $\rho$ admits a hard splitting into edge paths, the leftmost of which is $E_r \sigma_2$.

Suppose now that $H_q$ is a parabolic stratum.  Then arguing as in Lemma \ref{BasicHSplit}, we see that $\rho$ admits a hard splitting into edge paths, the leftmost of which is either $E_r \sigma_2$ or $E_r \sigma_2 \overline{E_q}$, where $\sigma_2$ has weight at most $q-1$.  Thus we may suppose that $\rho$ itself has this form.  Again, either the $D(n)$-step \nf \ of $\rho$ admits a hard splitting of the required form, or the $D(n)$-step nibbled future of $\rho$ has the form $E_r \sigma_3$, where $\sigma_3$ has weight at most $q-1$.  The above considerations cover the possibility that a \gep \ of weight $r$ occurs as a factor of the hard splitting.  Thus we may assume that in some \nf \ of $\rho$ there will necessarily be a hard splitting on each side of the edge of weight $r$.

In this fashion, going forwards into the \nf \ an amount of time bounded by a function of $n$, we may assume that $\rho$ has the form $E_r \sigma_4$, where $\sigma_4$ has weight exactly $s$ (if $\sigma_4$ has weight less than $s$ then $f_{\#}(E_r \sigma_4) = E_r \odot \sigma_5$ where $\sigma_5$ has weight less than $r$, and this is a splitting of the required form which is inherited by the \nf ).

We now consider what kind of stratum $H_s$ is.  Suppose that $H_s$ is parabolic.  There are only two ways in which cancellation between weight $s$ edges in the \nf \ of $\rho$ can occur (see \cite[Lemma 5.5]{BGroves}): they might be cancelled by edges whose immediate past is the edge of weight $r$ on the left end of the previous \nf ; alternatively, they can be nibbled from the right.  The $D(n)$-step entire future of $\rho$ admits a hard splitting as $E_r \odot \sigma_6$, where $\sigma_6$ has weight at most $r-1$.  There is no way that \nib \ can affect this splitting.

Finally, suppose that $H_s$ is an exponential stratum.  We follow a similar argument to the case when $H_r$ was an exponential stratum.  Either the $D(n)$-step \nf \ of $\rho$ admits a hard splitting of the desired kind (which means $\rho_{D(n)} = E_r \odot \sigma_7$ where $\sigma_7$ has weight at most $r-1$), or there are fewer $s$-illegal turns in the future of $\sigma_4$ in $\rho_{D(n)}$ than there are $s$-illegal turns in $\sigma_4$.  We then apply the same argument to the \nf \ of $\rho_{D(n)}$ until eventually we achieve a hard splitting of the required form.  This completes the proof of Proposition \ref{NibbleSplittingLemma}.
\end{proof}

We are now in a position to prove Theorem \ref{DecompTheorem}.  For this we require the following definition.

\begin{definition}
Suppose that $H_r$ is a stratum, and $E \in H_r$.
An {\em $r$-seed} is a non-empty subpath
$\rho$ of $f(E)$ which is maximal subject to lying in $G_{r-1}$.

If the stratum $H_r$ is not relevant, we just refer to {\em seeds}.
\end{definition}
Note that seeds are edge-paths and that the set of all seeds is finite.

The following is an immediate consequence of Lemma \ref{ExpHardSplitting} and RTT-(i) of Definition \ref{TrainTrackDef}.

\begin{lemma}
If $E \in H_r$ is an exponential edge and $\rho$ is an $r$-seed in
$f(E)$ then $f(E) = \sigma_1 \odot \rho \odot \sigma_2$ where
$\sigma_1$ and $\sigma_2$ are $r$-legal paths which start and finish
with edges in $H_r$.
\end{lemma}

\medskip

\begin{proof}[Proof (Theorem \ref{DecompTheorem}).]
Suppose that $\rho$ is a path of length $n$ and that $\rho_k$ is a
$k$-step \nf \ of $\rho$.  Denote by $\rho_0 = \rho, \rho_1, \dots ,
\rho_{k-1}$ the intermediate \nf s of $\rho$ used in order to define
$\rho_k$.  

We begin by constructing a van Kampen diagram\footnote{in fact, just a stack of corridors} $\Delta_k$ which encodes
the $\rho_i$, proceeding by induction
on $k$.  For $k = 1$ the diagram $\Delta_1$ has a single (folded)
corridor with the bottom labelled by $\rho$ and the path $\rho_1$ a
subpath of the top of this corridor.  Suppose that we have associated
a van Kampen diagram $\Delta_{k-1}$ to $\rho_{k-1}$, with a unique
corridor at each time $t = 0 , \ldots , k-2$, such that $\rho_{k-1}$
is a subpath of the top of the latest (folded) corridor.  Then we
attach a new folded corridor to $\Delta_{k-1}$ whose bottom is
labelled by $\rho_{k-1}$.  The path $\rho_k$ is, by definition, a
subpath of the top of this new latest corridor. By convention, we consider $\rho_i$ to occur at time $i$.

Choose an arbitrary edge $\e$ in $\rho_k$ on the (folded) top of the
latest corridor in $\Delta_k$.  We will prove that there is a path
$\sigma$ containing $\e$ in $\rho_k$ so that $\rho_k$ admits a hard
splitting immediately on either side of $\sigma$ and so that $\sigma$
is either suitably {\em short} or a \nf \ of a \gep .  The purpose of
this proof is to find a suitable notion of {\em short}. 

Consider the embedded `family forest' $\mathcal F$ for $\Delta_k$, tracing
the histories of edges lying on the folded tops of corridors (see \cite[3.2]{BGroves}).  
Let $p$ be the path in $\mathcal F$ which follows the history of
$\e$.  We denote by $p(i)$ the edge  which intersects $p$ and lies on
the bottom of the corridor at time $i$. The edges $p(i)$ form
the {\em past} of $\e$.  We will sometimes denote the edge $\e$ by
$p(k)$.  It will be an analysis of the times at which the weight of
$p(i)$ decreases that forms the core of the proof of the theorem. 

The weights of the edges $p(0), p(1), \ldots , p(k)$ form a
non-increasing sequence.  Suppose this sequence is $W = \{ w_0 ,
\ldots , w_k \}$. A {\em drop} in $W$ is a time $t$ such that $w_{t-1}
> w_t$.  At such times, the edge $p(t)$ is contained in a (folded)
seed in the bottom of a corridor of $\Delta_k$. 

We will show that either successive drops occur rapidly, or else we reach 
a situation wherein each time a drop occurs we lose no essential information by restricting our attention to a small subpath of $\rho_{i}$.

To make this localisation argument precise, we define
 {\em incidents}, which fall into two types.  
 
 An {\em incident of Type A} is a time $t$ which (i) is a drop; and (ii) is such that
there is a hard splitting of $\rho_t$ immediately on either side of
the folded seed containing $p(t)$.  

An {\em incident of Type B} is a
time $t$ such that $p(t-1)$ lies in an indivisible Nielsen path with a
hard splitting of $\rho_{t-1}$ immediately on either side, but $p(t)$
does not;  except that we do not consider this to be an incident if some $\rho_i$, for $i \le t-1$ admits a hard splitting $\rho_i = \sigma_1 \odot \sigma_2 \odot \sigma_3$ with $p(i) \subseteq \sigma_2$ and $\sigma_2$ a \gep .   In case of an incident of Type B, necessarily $p(t)$
lies in the nibbled future of a Nielsen path on one end of $\rho_t$ with a hard splitting of
$\rho_t$ immediately on the other side. 

Define the time $t_1$ to be the last time at which there is an
incident (of Type A or Type B).  If there are no incidents, let $t_1 =
0$.  If this incident is of Type A, the edge $p(t_1)$ lies in a folded 
seed, call it $\pi$, and there is a hard splitting of $\rho_{t_1}$
immediately on either side of $\pi$.  If the incident is of Type B,
the edge $p(t_1)$ lies in the $1$-step \nf \ of a Nielsen path, call
this \nf \ $\pi$ also.  In case $t_1 = 0$, let $\pi = \rho$.  We will
see that there is a bound, $\alpha$ say, on the length of $\pi$ which
depends only on $f$ and $n$, and not on the choice of $\pi$, or the
choice of \nf .  We postpone the proof of the existence of the bound
$\alpha$ while we examine the consequences of its existence. 

The purpose of isolating the path $\pi$ is that it is a path of
controlled length and the hard splitting \footnote{this splitting is vacuous in case
$t_1 = 0$ and at various other points during this proof which we do not
explicitly mention} of $\rho_{t_1}$
immediately on either side of $\pi$ means that we need only consider the \nf \ or
$\pi$.  Suppose that $\pi$ has weight $r$. 

{\bf Claim 1:}
There exists a constant $\beta = \beta(n,\alpha,f)$ so that one of the following must occur:  
\begin{enumerate}
\item[(i)] for some $t_1 < i < k$, the edge
$p(i)$ lies in  a \gep \ in $f_{\#}(\rho_{i-1})$
with a hard splitting immediately on either side; 
\item[(ii)] at some time $i \le t_1 + \beta$, the edge $p(i)$ 
lies in an indivisible Nielsen path $\tau$ in $f_{\#}(\rho_{i-1})$
with a hard splitting immediately on either side;  
\item[(iii)] $k - t_1 \le \beta$; or
\item[(iv)] there is a hard splitting of $\rho_k$ immediately on either side of $\e$. 
\end{enumerate}

\medskip

This claim implies the theorem, modulo the bound on $\alpha$, as we shall now explain.
In case (i), for all $j \ge i$, the edge $p(j)$ lies in the \nf \ of a \gep , so in particular this is true for $\e = p(k)$.  If case
(ii) arises then the definition of $t_1$ implies that for $j \ge i$, the edge $p(j)$ always
lies in a path labelled $\tau$ with a hard splitting immediately on
either side, for otherwise there would be a subsequent incident.  Also, the length of this Nielsen path is at most $\alpha L^\beta$.  If case (iii) arises, then the nibbled future of
$\pi$ at time $k$ has length at most $\alpha L^{\beta}$.

To prove the claim, we define two sequences of numbers
$V_{\omega}, V_{\omega - 1}, \ldots , V_1$ and $V'_{\omega},
V'_{\omega - 1}, \ldots , V'_1$, depending on $n$ and $f$, as follows (where $D'(n)$ is the function
  from Proposition \ref{NibbleSplittingLemma}): 
\begin{eqnarray*}
V_{\omega} &  := & D'(\alpha),\\ 
V'_{\omega} & := & V_{\omega} + \alpha L^{V_{\omega}}.
\end{eqnarray*}
For $\omega > i \ge 1$, supposing $V'_{i+1}$ to be defined, 
\[      V_i   :=   V'_{i+1} + D'(\alpha L^{V'_{i+1}}).  \] 
Also, supposing  $V_i$ to be defined, we define
\[      V'_i   :=  V_i +  \alpha L^{V_i} .      \]

We consider the situation at time $t_1 + V_r$ (recall that $r$ is the
weight of $\pi$).  Possibly $k \le t_1 + V_r$, which is covered by
case (iii) of our claim.  Therefore, suppose that $k > t_1 + V_r$. 

According to Proposition \ref{NibbleSplittingLemma}, and the definition of
$t_1$,  at time $t_1 + V_r$ the $V_r$-step \nf \ of $\pi$ which exists
in $\rho_{t_1 + V_r}$ admits a hard splitting into edge paths, each of
which is either: 
\begin{enumerate}
\item a single edge of weight $r$;
\item   a \nf \ of a weight $r$ indivisible Nielsen path;
\item   a \nf \ of a weight $r$ \gep ; or
\item   a path of weight at most $r-1$.
\end{enumerate}
We need to augment possibility (3) by noting that the proof of Proposition \ref{NibbleSplittingLemma} shows that the \gep \ referred to lies in the $j$-step \nf \ of $\pi$ for some $j \le V_r$.

We analyse what happens when the edge $p(t_1+V_r)$ lies in each of these four types of path.

{\bf Case (1):}
In the first case, by the definition of $t_1$, there will be
a hard splitting of $\rho_k$ immediately on either side of $\e$, since
in this case if there is a drop in $W$ after $t_1 + V_r$ then there is an
incident of Type A, contrary to hypothesis. 

{\bf Case (3):}
If $p(t_1+V_r)$ lies in a
path of the third type then we are in case (i) of our claim, and hence
content.  

The fourth type of path will lead us to an inductive
argument on the weight of the path under consideration.  But first we
consider the \nf s of Nielsen paths. 

{\bf Case (2):}
Suppose that in $\rho_{t_1 + V_r}$ the edge $p(t_1 + V_r)$ lies in the
\nf \ of a Nielsen path of weight $r$, with a hard splitting of
$\rho_{t_1 + V_r}$ immediately on either side.  Suppose that this \nf
\ is $\pi_r$.  If $\pi_r$ is actually a Nielsen path then we lie in
case (ii) of our claim.  Thus suppose that $\pi_r$ is not a Nielsen
path.  It has length at most $\alpha L^{V_r}$, and within  time
$\alpha L^{V_r}$ any \nf \ of $\pi_r$ admits a hard splitting into edge paths of
types (1), (3) and (4) from the above list.  

To see this, consider the
three types of indivisible Nielsen paths.  If $\tau$ is a Nielsen path
which is a single edge fixed pointwise by $f$, then any \nf \ of
$\tau$ is either a single edge or empty.  Suppose that $\tau$ is an
indivisible Nielsen path of weight $r$ and $H_r$ is exponential, and
suppose that $\tau'$ is a proper subpath of $\tau$.  Then there is
some iterated image $f_{\#}^l(\tau')$ of $\tau'$ which is $r$-legal with $l < \alpha L^{V_r}$.
Finally suppose that $E_i \tau^k \overline{E_j}$ is an indivisible
Nielsen path of parabolic weight.  Thus $\tau$ is a Nielsen path of
weight less than $r$, and $E_i$, $E_j$ are edges such that $f(E_i) =
E_i \odot \tau^m$, $f(E_j) = E_j \odot \tau^m$.  A $1$-step \nf \ of
$E_i \tau^k \overline{E_j}$ has one of three forms:  (I) $E_i
\tau^{k_1} \tau'$, where $\tau'$ is a proper sub edge-path of $\tau$;
(II) $\tau' \tau^{k_2} \tau''$ where $\tau'$ and $\tau''$ are proper
sub edge-paths of $\tau$; or (III) $\tau' \tau^{k_3} \overline{E_j}$,
where $\tau'$ is a proper sub edge-path of $\tau$. 

\smallskip

{\em Case 2(I):}
In this case, $E_i \tau^{k_1} \tau'$ admits a hard splitting into $E_i$ and $\tau^{k_1}
\tau'$, which is of the required sort.  

\smallskip

{\em Case 2(II):}
In this case the path already
had weight less than $r$.  

\smallskip

{\em Case 2(III):}
Suppose we are in case (III), and that $\mu$, the $\alpha L^{V_r}$-step \nf \ of $\tau'\tau^{k_3} \overline{E_j}$  has a copy of $\overline{E_j}$. Lemma \ref{NibNoTau} assures us that no \nf \ of $\tau'$ can contain $\tau$ as a subpath, and therefore there is a splitting of $\mu$ immediately on the right of $\overline{E_r}$, and we are done.  If there is no copy of $\overline{E_j}$ in $\mu$, we are also done, since this \nf \ must have weight less than $r$.

{\bf Case 4:}
Having dealt with cases (1) and (3), we may now suppose that at time $t_1 + V_r + \alpha L^{V_r} = t_1 +
V'_r$ the edge $p(t_1 + V'_r)$ lies in an edge-path of weight at most
$r-1$ with a hard splitting of $\rho_{t_1 + V'_r}$ immediately on
either side.\footnote{Note that again it is possible that $k < t_1 +
  V'_r$, in which case we are in case (iii) of our claim.  We suppose
  therefore that this is not the case.}  Denote this path by $\pi'_r$, chosen to be in the future of $\pi$.
Note that $\pi'_r$ has length at most $\alpha L^{V'_r}$. 

By Proposition \ref{NibbleSplittingLemma} again, either $k < t_1 + V_{r-1}$
or at time $t_1 + V_{r-1}$ the \nf \ of $\pi'_r$ admits a hard
splitting into edge-paths each of which is either: 
\begin{enumerate}
\item   a single edge of weight $r-1$;
\item   a \nf \ of a weight $r-1$ indivisible Nielsen path;
\item   a \nf \ of a weight $r-1$ \gep ; or
\item   a path of weight at most $r-2$.
\end{enumerate}

We continue in this manner.  We may conceivably fall into case (4)
each time until $t_1 + V_1$ when it is not possible to fall into a
path of weight at most $1 - 1$!  Thus at some stage we must fall into
one of the first three cases.  This completes the proof of Claim 1.

\bigskip 

{\bf The existence of $\alpha$}.
We must find a
bound, in terms of $n$ and $f$, on the length of indivisible Nielsen
paths that arise in the \nf \ of $\rho$ with a hard splitting
immediately on either side.\footnote{Recall that the definition of
  Type B incidents excluded the case of Nielsen paths which lie in the
  \nf \ of a \gep \ with a hard splitting immediately on either side.}
To this end, suppose that $\e'$ is an edge which lies in an
indivisible Nielsen path $\tau$ in a $k'$-step \nf \ of $\rho$, and
that there is a hard splitting immediately on either side of $\tau$. 
We again denote the $i$-step nibbled future of $\rho$ by $\rho_i$ for
$0 \le i \le k'$. 

As above, we associate a diagram $\Delta_{k'}$ to
$\rho_{k'}$.\footnote{If we are considering Nielsen paths arising in
  the past of $\e$ above, then we can assume $k' \le k$ and that
  $\Delta_{k'}$ is a subdiagram of $\Delta_k$ in the obvious way.}
Denote by $q$ the path in the family forest of $\Delta_{k'}$ which follows the past of $\e'$. Let $q(i)$ be the edge in $\rho_i$ which intersects $q$.  Let the sequence
of weights of the edges $q(i)$ be $W' = \{ w'_0, \ldots , w'_{k'} \}$. 

Define incidents of Type A and B for $W'$ in exactly the same way as
for $W$, and let $t_2$ be the time of the last incident of Type A for
$W'$.  If there is no incident of Type A for $W'$ let $t_2 = 0$.  Let
$\kappa$ be the folded seed containing $q(t_2)$; in case $t_2 = 0$
let $\kappa = \rho$.  Define $\theta = \mbox{max} \{ n , L \}$ and note
that $| \kappa | \le \theta$.  The path $\tau$ must lie in the \nf \ of
$\kappa$, so it suffices to consider the \nf \ of $\kappa$.  Suppose
that $\kappa$ has weight $r'$. 

We deal with the \nf \ of $\kappa$ in the same way as we dealt with
that of $\pi$.   Let $\kappa_0 = \kappa, \kappa_1, \ldots$ be the \nf s of $\kappa$.

\medskip

{\bf Claim 2:} There exists a constant $\beta' = \beta'(n,f)$ so that one of the following must occur:  
\begin{enumerate}
\item[(i)] for some $t_2 < i < k'$, the edge $q(i)$ lies in a \gep \ in $f_{\#}(\kappa_{i-1})$ that has a hard splitting immediately on either side; 
\item[(ii)] not in case (i), {\em and} at some time $i \le k'$ the edge $q(i)$ lies in an indivisible Nielsen
path $\tau_0$ in $f_{\#}(\kappa_{i-1})$ so that $|\tau_0| \le\theta L^{\beta'}$ and immediately on either side of $\tau_0$ there is a hard splitting, {\em and} there are no incidents of Type B after time $i$;
\item[(iii)] $k' - t_2 \le \beta'$; or 
\item[(iv)] there is a hard splitting of $\kappa_{k'}$ immediately on either side of $\e'$. 
\end{enumerate}

Let us prove that this claim implies the existence of $\alpha$ and hence completes the proof of the theorem.  By definition, $\alpha$ is required to be an upper bound on the length of an arbitrary Nielsen path $\tau$ involved in a Type B incident.  We assume this incident occurs at time $k'$ and use Claim 2 to analyse what happens.  

Case (i) of Claim 2 is irrelevant in this regard.  If case (ii) occurs, the futures of $\tau_0$ are unchanging up to time $k'$, so $\tau = \tau_0$ and we have our required bound.  In case (iii) the length of $\tau$ is at most $\theta L^{\beta'}$, and in case (iv) $\tau$ is a single edge.  It suffices to let $\alpha = \theta L^{\beta'}$.

\medskip

It remains to prove Claim 2.  The proof of Claim 2 follows that of Claim 1 almost verbatim, with $\theta$ in place of $\alpha$ and $\kappa$ in place of $\rho$, etc.,  {\em except}  that the third sentence in Case (2) of the proof becomes invalid because Type B incidents after time $t_2 + V_r$ may occur.  

In this setting, suppose $\pi_r$ (which occurs at time $t_2 + V_r$) is a Nielsen path, but that we are not in case (ii) of Claim 2, and there is a subsequent Type B incident at time $j$, say.  The length of $\pi_r$ is at most $\theta L^{V_r}$. The Nielsen path at time $j-1$ has the same length as the one at time $t_2 + V_r$.  We go forward to time $j$, where the future of $\pi_r$ is no longer a Nielsen path, and continue the proof of Case (2) from the fourth sentence of the proof.

Otherwise, the proof of Claim 2 is the same as that of Claim 1 (the above modification is required at each weight, but at most {\em once} for each weight).  The only way in which the length bounds change is in the replacement of $\theta$ by $\alpha$ (including in the definitions of $V_i$ and $V'_i$). This finally completes the proof of Theorem \ref{DecompTheorem}.
\end{proof}

\section{Passing to an iterate of $f$} \label{IterateSection}

In this section we describe what happens to various definitions when we replace 
$f$ by an iterate.  Suppose that $k \ge 1$, and consider the relationship between 
$f$ and $f_0 = f_{\#}^k$.

First, for any integer $j \ge 1$, the set of $kj$-monochromatic paths for $f$ is the same
as the set of $j$-monochromatic paths for $f_0$.  Therefore, once Theorem 
\ref{MainProp} is proved, we will pass to an iterate so that $r$-monochromatic
becomes $1$-monochromatic. However, the story is not quite as simple as that.

It is not hard to see that if $\sigma \odot \nu$ is a 
hard splitting for $f$, then it is also a hard splitting for $f_0$.

When $f$ is replaced by $f_0$, the set of \gep s is unchanged, as are the sets of 
\pep s and indivisible Nielsen paths.  Also, the set of indivisible
Nielsen paths which occur as sub-paths of $f(E)$ for some linear edge $E$ remains unchanged.

With the definition as given, the set of $(J,f_0)$-atoms may be smaller than the set of $(J,f)$-atoms.  This is because an atom is required to be 
$1$-monochromatic.  However, we continue to consider the set of $(J,f)$-atoms
even when we pass to $f_0$, and we also consider paths to be {\em beaded}
if they are $(J,f)$-beaded.

Since we are quantifying over a smaller set of paths the constant $V(n,f_0)$ in Theorem \ref{DecompTheorem} is assumed, without loss of generality, to be $V(n,f)$.  This is an 
important point, because the constant $V$ is used to find the appropriate $J$ when proving 
Theorem \ref{MainProp}.  When passing from $f$ to $f_0$, we need this $J$ to 
remain unchanged, for the appropriate iterate $k$ which we eventually choose depends 
crucially upon $J$.  See Remark \ref{SameJ} below.

It is also clear that if $m \le n$ then without loss of generality we may assume that $V(m,f) \le V(n,f)$.  Once again, this is because we are considering a smaller set of paths when defining $V(m,f)$.

We now want to replace $f$ by a fixed iterate in order to control some of the cancellation within monochromatic paths.  The following lemma is particularly useful in the proof of Proposition \ref{gammasingleedge} below, and also for Theorem 
\ref{ColourCancelMain}.

\begin{lemma} \label{Power1}
There exists $k_1 \ge 1$ so that $f_1 = f_{\#}^{k_1}$ satisfies the following. 
Suppose that $E$ is an exponential edge of weight $r$ and that $\sigma$ is an 
indivisible Nielsen path of weight $r$ (if it exists, $\sigma$ is unique up to a change 
of orientation). Then
\begin{enumerate}
\item $|f_1(E)| > |\sigma|$;
\item Moreover, if $\sigma$ is an indivisible Nielsen path of exponential weight 
$r$ and $\sigma_0$ is a proper subedge-path of $\sigma$,  then $(f_1)_{\#}(\sigma_0)$ is $r$-legal;
\item If $\sigma_0$ is a proper initial sub edge-path of $\sigma$ then $(f_1)_{\#}(\sigma_0)$ admits a hard splitting, $f(E) \odot \xi$, where $E$ is the edge on the left end of $\sigma$;
\item Finally, if $\sigma_1$ is a proper terminal sub edge-path of $\sigma$ then $(f_1)_{\#}(\sigma_1) = \xi' \odot f(E')$ where $E'$ is the edge on the right end of $\sigma$.
\end{enumerate}
Now suppose that $\sigma$ is an indivisible Nielsen path of parabolic weight $r$ and that $\sigma$ is a sub edge-path of $f(E_1)$ for some linear edge $E_1$.  The path $\sigma$ is either of the form $E \eta^{m_{\sigma}} \overline{E'}$ or of the form $E \overline{\eta}^{m_{\sigma}} \overline{E'}$, for some linear edges $E$ and $E'$.  Then 
\begin{enumerate}
\item If $\sigma_0$ is a proper initial sub edge-path of $\sigma$ then
\[      (f_1)_{\#}(\sigma_0) = E \odot \eta \odot \cdots \odot \eta \odot \xi''   ,       \]
where there are more than $m_{\sigma}$ copies of $\eta$ visible in this splitting.
\item If $\sigma_1$ is a proper terminal sub edge-path of $\sigma$ then
\[      (f_1)_{\#}(\sigma_1) = \xi' \odot \overline{\eta} \odot \cdots \odot \overline{\eta} \odot \overline{E'},    \]
where there are more than $m_\sigma$ copies of $\overline{\eta}$ visible in this splitting;
\end{enumerate}
\end{lemma}
\begin{proof}
First suppose that  $H_r$ is an exponential stratum, that $\sigma$ is an indivisible Nielsen path of weight $r$, and that $E$ is an edge of weight $r$.  Since $|f_{\#}^j(E)|$ grows exponentially with $j$, and $|f_{\#}^j(\sigma)|$ is constant, there is certainly some $d_0$ so that $|f_{\#}^d(E)| > |\sigma|$ for all $d \ge d_0$.

There is a single $r$-illegal turn in $\sigma$, and if $\sigma_0$ is a proper sub edge-path of $\sigma$.  By Lemma \ref{NoTau}, no future of $\sigma_0$ can contain $\sigma$ as a subpath.  The number of $r$-illegal turns in iterates of $\sigma_0$ must stabilise, so by Lemma \ref{ExpSplitting} there is an iterate of $\sigma_0$ which is $r$-legal.  Since there are only finitely many paths $\sigma_0$, we can choose an iterate of $f$ which works for all such $\sigma_0$.

Suppose now that $\sigma_0$ is a proper initial sub edge-path of $\sigma$, and that $E$ is the edge on the left end of $\sigma$.  It is not hard to see that every (entire) future of $\sigma_0$ has $E$ on its left end.  We have found an iterate of $f$ so that $f_{\#}^{d'}(\sigma_0)$ is $r$-legal.  It now follows immediately that
\[      f_{\#}^{d'+1}(\sigma_0) = f(E) \odot \xi        ,       \]
for some path $\xi$.  The case when $\sigma_1$ is a proper terminal sub edge-path of $\sigma$ is identical.

Now suppose that $H_r$ is a parabolic stratum and that $\sigma$ is an indivisible Nielsen path of weight $r$ of the form in the statement of the lemma.  The claims about sub-paths of $\sigma$ follow from the hard splittings $f(E) = E \odot u_E$ and $f(E') = E' \odot u_{E'}$, and from the fact that $m_{\sigma}$ is bounded because $\sigma$ is a subpath of some $f(E_1)$.

As in Remark \ref{stable}, we can treat each of the cases separately, and finally take a maximum.  
\end{proof}

\section{The \nf s of \gep s} \label{nfGEPSection}

The entire future of a \gep \ is a \gep \ but a \nf \ need not be and Theorem \ref{DecompTheorem} tells us that we need to analyse these \nf s.  This analysis will lead us to define {\em proto-\pep s}.  In Proposition \ref{gammasingleedge}, we establish a
normal form for proto-\pep s which proves that proto-\pep s are in fact the \pep s which 
appear in the Beaded Decomposition Theorem.

To this end, suppose that
\[ \zeta = E_i \overline\tau^n \overline{E_j}           \]
is a \gep , where $\tau$ is a Nielsen path, $f(E_i) = 
E_i \odot \tau^{m_i}$ and $f(E_j) = E_j \odot \tau^{m_j}$.  As in Definition
\ref{Exceptional}, we consider $E_i \overline\tau^n \overline{E_j}$ to be
unoriented, but here we do not suppose that $j \le i$.  However, we suppose $n > 0$ and thus,
since $E_i \overline\tau^n \overline{E_j}$ is a \gep , $m_j > m_i > 0$.

The analysis of \gep s of the form $E_j {\tau}^n \overline{E_i}$ is entirely similar
to that of \gep s of the form $E_i \overline\tau^n \overline{E_j}$ except that one must reverse all left-right orientations.  Therefore, we ignore this case until Definition \ref{pep} below (and often afterwards also!).

We fix a sequence of \nf s $\zeta = \rho_{-l}, \ldots , \rho_0, \rho_1, \ldots , \rho_k, \ldots$ of $\zeta$, where $\rho_0$ is the first \nf \ which is not the entire future.  Since the entire future of a \gep \ is a \gep , we restrict our attention to the \nf s of $\rho_0$.

There are three cases to consider, depending on the type of sub-path on either end of $\rho_0$.
\begin{enumerate}
\item $\rho_0 = \bar{\tau_0} \bar{\tau}^m \overline{E_j}$;
\item $\rho_0 = \bar{\tau_0} \bar{\tau}^m \bar{\tau_1}$.
\item $\rho_0 = E_i \bar{\tau}^m \bar{\tau_1}$;
\end{enumerate}
where $\tau_0$ is a (possibly empty) initial sub edge-path of $\tau$, and $\tau_1$ is
a (possibly empty) terminal sub-edge path of $\tau$.

In case (1) $\rho_0$ admits a hard splitting
\[      \rho_0 = \overline\tau_0 \odot \overline\tau \odot \cdots \odot \overline\tau \odot \overline E_j.    \]
Since  $\tau_0$ is a sub edge-path of $f(E_i)$, it has length less than
 $L$ and its nibbled futures admit hard splittings as in
Theorem \ref{DecompTheorem} into \nf s of GEPs and paths of length at most
$V(L,f)$. These \gep s 
 will necessarily be of strictly lower weight than $\rho_0$, since $\overline\tau_0$ is. Thus, 
case (1) is easily dealt with  by an induction on weight, supposing that we have a nice 
splitting of the \nf s of lower weight \gep s; this is made precise in Proposition \ref{nfgeptopep}.  Case (2) is entirely similar.

Case (3)  is by far the most troublesome of the three, and it is this case which leads to the definition of {\em proto-\pep s} in Definition \ref{pep} below.  Henceforth assume $\rho_0 = E_i \bar{\tau}^m \bar{\tau_1}$.

Each of the \nf s of $\rho_0$ (up to the moment of death, Subsection \ref{Death}) 
has a \nf \ of 
$\overline\tau_1$ on the right.  If the latter becomes empty at some point, 
the \nf \ of $\rho_0$ at this time has the form $E_i \overline\tau^{n'} \overline\tau_2$, where $\tau_2$ is a proper sub edge-path of  $\tau$.  We
 restart our analysis at this moment.  Hence we make the following

\begin{workingassumption} \label{NibAssumption}
We make the following two assumptions on the \nf s considered:
\begin{enumerate}
\item All nibbling of $\rho_i$ occurs on the right; and 
\item the $i$-step \nf \ $\overline\tau_{1,i}$ of $\overline\tau_1$ inherited from $\rho_i$ is non-empty.
\end{enumerate}
\end{workingassumption}

We will deal with the case $m - km_i < 0$ later, in particular with the value of $k$ for 
which $m - (k-1)m_i \ge 0$ but $m - km_i < 0$.  For now suppose that $m - km_i \ge 0$.

In this case, the path $\rho_k$ 
has the form
\[      \rho_k = E_i \overline\tau^{m - km_i} \overline\tau_{1,k} .     \]

There are (possibly empty) Nielsen edge paths $\iota$ and $\nu$, and an indivisible Nielsen edge path $\sigma$ so that 
\begin{equation} \label{TauDecomp}
\tau = \iota \odot \sigma \odot \nu \mbox{ and } \tau_1 = \sigma_1 \odot \nu     ,
\end{equation}
where $\sigma_1$ is a proper terminal sub edge-path of $\sigma$.  Now, as in Working Assumption \ref{NibAssumption}, there is no loss of generality in supposing that
\[      \rho_k = E_i \overline\tau^{m-km_i} \bar \nu \bar\sigma_{1,k} , \]
where $\overline\sigma_{1,k}$ is the \nf \ of $\overline\sigma_1$ inherited from $\rho_k$, and that $\overline\sigma_{1,k}$ is non-empty.

Since $|\sigma_1| < L$, by Theorem \ref{DecompTheorem} the path $\sigma_{1,k}$ admits a hard splitting into edge-paths each of which is either the \nf \ of a \gep , or of length at most $V(L,f)$;  we take the (unique) maximal hard splitting of $\sigma_{1,k}$ into edge paths.  

Let $s = \lfloor m / m_i \rfloor + 1$.  In $\rho_s$ (but not before) there may be some 
interaction between the future of $E_i$ and $\overline\sigma_{1,s}$.   We denote by $\gamma_{\sigma_1}^{k,m}$ the concatenation of those factors in the hard splitting of $\overline\sigma_{1,k}$
which contain edges any part of whose future is eventually cancelled by some edge in the future 
of $E_i$ under any choice of \nf s of $\rho_k$ (not just the $\rho_{k+t}$ chosen earlier) and any choice of tightening.  Below we will analyse more carefully the structure of the paths $\overline\sigma_{1,k}$ and $\gamma_{\sigma_1}^{k,m}$.

We now have $\overline\sigma_{1,k} = \gamma_{\sigma_1}^{k,m} \odot \sigma_{1,k}^\bullet$.  From (\ref{TauDecomp}), we also have
\begin{equation} \label{nfpep}
\rho_k = E_i \overline\tau^{m-km_i}  \overline\nu \gamma_{\sigma_1}^{k,m} \odot \sigma_{1,k}^\bullet .
\end{equation}

\begin{definition}[Proto-\pep s] \label{pep}
Suppose that $\tau$ is a Nielsen edge path, $E_i$ a linear edge such that $f(E_i) = E_i 
\odot \tau^{m_i}$ and $\tau_1$ a proper terminal sub edge-path of $\tau$ such that $\tau_1 
= \sigma_1 \odot \nu$ as in (\ref{TauDecomp}).  Let $k,m \ge 0$ be such that $m - km_i \ge 0$ and let $\gamma_{\sigma_1}^{k,m}$ be as in (\ref{nfpep}).
A path $\pi$ is called a {\em proto-\pep} if either $\pi$ of $\overline{\pi}$ is of the form
\[      E_i  \overline\tau^{m-km_i} \overline\nu \gamma_{\sigma_1}^{k,m} .       \]
\end{definition}

\begin{remarks} 
\ \par

\begin{enumerate}
\item The definition of proto-\pep s is intended to capture those paths which remain when a \gep \ is 
partially cancelled, leaving a path which may shrink in size of its own accord.
\item By definition, a proto-\pep \ admits no non-vacuous hard splitting into edge paths.
\end{enumerate}
\end{remarks}

We now introduce two distinguished kinds of proto-\pep s.

\begin{definition} \label{stabletransient}
Suppose that 
\[      \pi = E_i  \overline\tau^{m-km_i} \overline\nu \gamma_{\sigma_1}^{k,m} ,\] 
is a proto-\pep \ as in Definition \ref{pep}.

The path $\pi$ is a {\em transient} proto-\pep \  if $k=0$.

The path $\pi$ is a {\em stable} proto \pep \ if $\gamma_{\sigma_1}^{k,m}$
is a single edge.
\end{definition}

\begin{lemma} \label{transientispep}
A transient proto-\pep \ is a \pep. 
\end{lemma}
\begin{proof}
With the notation of Definition \ref{pep}, in this case $\gamma^{0,m}_{\sigma_1}$
is visibly a sub-path of $\bar{\tau}$, and the proto-\pep \ is visibly a sub-path of a \gep.
\end{proof}

\begin{lemma} \label{stableprotoispep}
A stable proto-\pep \ is a \pep.
\end{lemma}
\begin{proof}
Since $\bar{\sigma}$ is a Nielsen path, if $\alpha$ is a \nf \ of $\bar{\sigma}$ where
all the nibbling has occurred on the right, then the first edge in $\alpha$ is the
same as the first edge in $\bar{\sigma}$.

However, $\gamma_{\sigma_1}^{k,m}$ is a \nf \ of $\bar{\sigma}$ where all the
nibbling has occurred on the right.  Therefore, if $\gamma_{\sigma_1}^{k,m}$
is a single edge then it must be a sub-path of $\bar{\sigma}$.  It is now easy to
see that any stable proto-\pep \ must be a \pep.
\end{proof}

\begin{remark}
We will prove in Proposition \ref{gammasingleedge} that after
replacing $f$ by a suitable iterate all proto-\pep s are either
transient or stable, and hence are \pep s.
\end{remark}

\subsection{The Death of a proto-\pep } \label{Death}

Suppose that $\pi = E_i \overline\tau^{m-km_i} \overline\nu \gamma_{\sigma_1}^{k,m}$
 is a proto-\pep \ which satisfies Assumption \ref{NibAssumption}.  Let $q = \lfloor \frac{m-km_i}{m_i}
 \rfloor + 1$, and consider, $\pi_{q-1}$, a $(q-1)$-step \nf \  of $\pi$. As before, we assume 
 that the $(q-1)$-step \nf \ of $\gamma_{\sigma_0}^{k,m}$ inherited from a $\pi_{q-1}$
  is not empty and that the edge labelled $E_i$ on the very left is not nibbled. 

In $\pi_{q-1}$, the edge $E_j$ has consumed all of the copies of $\overline\tau$ and begins to interact 
with the future of $\overline\nu \gamma_{\sigma_1}^{k,m}$. Also, the future of $\pi$ at time $q$ need not contain a \pep.  Hence we refer to the time $q$ as the {\em death of the \pep}.
Recall that $\tau = \iota \odot \sigma 
\odot \nu$ and that $\gamma_{\sigma_1}^{k,m}$ is a $k$-step \nf \ of $\overline\sigma_1$, where $\sigma_1$ is a proper subbath of $\sigma$.  Let $p = m - (k+q-1)m_i $, so that $0 \le p < m_i$.

The path $\pi_{q-1}$ has the form
\[      \pi_{q-1} = E_i \overline\tau^{p} \overline\nu \gamma_{\sigma_1}^{k+q-1,m} .       \]
Suppose that $\pi_q$ is a $1$-step
\nf \ of $\pi_{q-1}$.  In other words, $\pi_q$ is a subpath of $f_{\#}(\pi_{q-1})$.
Consider what happens when $f(\pi_{q-1})$ is tightened to form $f_{\#}(\pi_{q-1})$ 
(with any choice of tightening).   The $p$ copies of $\overline\tau$ 
(possibly in various stages of tightening) will be consumed by $E_i$, leaving 
$\overline \nu \odot f(\gamma_{\sigma_1}^{k+q-1,m})$ to interact with at least one 
remaining copy of $\tau = \iota \odot \sigma \odot \nu$.  The paths $\nu$ and 
$\overline{\nu}$ will cancel with each other\footnote{The hard splittings imply that 
this cancellation must occur under {\em any} choice of tightening.}.

Lemma \ref{NibNoTau} states that $\gamma_{\sigma_1}^{k,m}$ cannot contain $\sigma$ as a subpath.  Therefore, once $\nu$ and $\overline{\nu}$ have cancelled, not all of $\overline{\sigma}$ will cancel with $f(\gamma_{\sigma_1}^{k+q-1,m})$.  A consequence of this discussion (and the fact that $f(E_i) = E_i \odot \tau^{m_i}$) is the following

\begin{lemma} \label{pepSplitting}
Suppose that $\pi = E_i \overline\tau^{m-km_i} \overline\nu \gamma_{\sigma_0}^{k,m}$ is a proto-\pep , and let $q = \lfloor \frac{m-km_i}{m_i} \rfloor + 1$.  Suppose that $\pi_{q-1}$ is a $(q-1)$-step \nf \ of $\pi$ satisfying Assumption \ref{NibAssumption}.  If $\pi_q$ is an immediate \nf \ of $\pi_{q-1}$ and $\pi_q$ contains $E_i$ then $\pi_q$ admits a hard splitting
\[      \pi_q = E_i \odot \lambda.      \]
\end{lemma}

We now analyse the interaction between $f(\gamma_{\sigma_1}^{k+q-1,m})$ and $\sigma$ more closely.  As usual, there are two cases to consider, depending on whether $\sigma$ has exponential or parabolic weight\footnote{Recall that there are three kinds of indivisible Nielsen paths: constant edges, parabolic weight and exponential weight.  If $\sigma$ has nontrivial proper sub edge-paths, then it is certainly not a single edge, constant or not.}.

In the following proposition, $f_1$ is the iterate of $f$ from Lemma \ref{Power1} and we are using the definitions as explained in Section \ref{IterateSection}.  Also, we assume that proto-\pep s are defined using $f_1$, not $f$.

\smallskip

\begin{proposition} \label{gammasingleedge}
Every proto-\pep \ for $f_1$ is either transient or stable.
In particular, every proto-\pep \  for $f_1$ is a \pep.
\end{proposition}
\begin{proof}
Let $\pi = E_i \overline\tau^{m-km_i} \overline\nu \gamma_{\sigma_1}^{k,m}$ be a 
proto-\pep \ for $f_1$.

Lemma \ref{transientispep} implies that  if $k = 0$ then $\pi$ is a \pep.
Consider Working Assumption \ref{NibAssumption}.
  If Assumption \ref{NibAssumption}.(2) fails to hold at any point, then we can restart our analysis, and 
  in particular we have a transient proto-\pep \  at this moment.
  Thus we may suppose that $\pi$ is an initial sub-path of a $k$-step \nf \ of a \gep , where 
  $k \ge 1$ and we may further suppose that $\pi$ satisfies Assumption \ref{NibAssumption}.(2). We prove that in this case $\pi$ is a stable proto-\pep .

First suppose that $\sigma$ has exponential weight, $r$ say.  If $\overline{\sigma}_0$ is a 
proper initial sub edge-path of $\overline\sigma$ then Lemma \ref{Power1} asserts that 
\[      (f_1)_{\#}(\sigma_0) = f(E) \odot \xi   ,       \]
and $| f(E)| > |\sigma|$.  Note also that $f(E) = E \odot \xi''$ for some path $\xi''$.

Now, at the death of the proto-\pep , the \nf \ of $\gamma_{\sigma_0}^{k,m}$ interacts 
with a copy of $E_i$, and in particular with a copy of $f(\sigma)$ (in some stage
 of tightening).  Now the above hard splitting, and the fact that $\sigma$ is not 
 $r$-legal whilst $f(E)$ is, shows that $\gamma_{\sigma_1}^{k,m}$ must be a 
 single edge (namely $E$).

Suppose now that $\sigma$ has parabolic weight.  Since $\sigma$ has proper 
sub edge-paths, it is not a single edge and so $\sigma$ or $\overline{\sigma}$ 
has the form $E \eta^{m_{\sigma}} \overline{E'}$.  The hard splittings guaranteed
 by Lemma \ref{Power1} now imply that $\gamma_{\sigma_1}^{k,m}$ is a single
  edge in this case also.

Therefore, every proto-\pep \ for $f_1$ is transient or stable, proving the first
assertion of the proposition.  The second assertion follows from the first assertion,
and Lemmas \ref{transientispep} and \ref{stableprotoispep}.
\end{proof}

Finally, we can prove the main result of this section.  In the following, $L_1$ is the maximum length of $f_1(E)$ over all edges $E$ of $G$.

The following statement assumes the conventions of Section \ref{IterateSection}.
\begin{proposition} \label{nfgeptopep}
Under iteration of the map $f_1$ constructed in Lemma \ref{Power1}, any
\nf \ of a \gep \ admits a hard splitting into edge paths, each of which is either a 
\gep , a \pep , or of length at most $V(2L_1,f)$.
\end{proposition}
\begin{proof}
Suppose that $E_i \overline{\tau}^n \overline{E_j}$ is a \gep \ of weight $r$.  We may suppose
 by induction that any \nf \ of any \gep \ of weight less than $r$ admits a hard splitting of
the required form (the base case $r=1$ is vacuous, since there cannot be a \gep \ of weight 
$1$).

Suppose that $\rho$ is a \nf \ of $E_i \tau^n \overline{E_j}$.  If $\rho$ is the entire future,
it is a \gep \ and there is nothing to prove.  Otherwise, as in the analysis at the beginning of this section, 
we consider the first time when a \nf \ is not the entire future.  Let the \nf \ be $\rho_0$.  
In cases (1) and (2) from that anaylsis, $\rho_0$ admits a hard splitting into edge-paths, 
each of which is either (i) $\overline{E}_i$; (ii) $\overline{\tau}$; or (iii) a proper 
sub edge-path of $\tau$.  In each of these cases, Theorem \ref{DecompTheorem} asserts
 that there is a hard splitting of $\rho$ into edge-paths, each of which is either of length at
  most $V(L,f)$ or is the \nf \ of a \gep .  Any \nf \ of a \gep \ which occurs in this splitting is 
  necessarily of weight strictly less than $r$, and so admits a hard splitting of the required 
  form by induction.  

Suppose then that $\rho_0$ satisfies Case (3), the third of the cases articulated at 
the beginning of this section.  In this 
case, $\rho_0$ is a transient proto-\pep .  Also, any time that Assumption \ref{NibAssumption}.(2) 
is not satisfied, the \nf \ of $\rho_0$ is a transient proto-\pep .  Thus, we may assume that 
Assumption \ref{NibAssumption} is satisfied. If $m -km_i \ge 0$ then we have
\[      \rho = E_i \overline\tau^{m-km_i} \overline\nu \gamma_{\sigma_1}^{k,m} \odot 
\sigma_{1,k}^\bullet    .       \]
The first path in this splitting is a stable \pep \ by Proposition \ref{gammasingleedge}.  
Once again, Theorem \ref{DecompTheorem} and the inductive hypothesis yield a hard 
splitting of $\sigma_{1,k}^\bullet$ of the required form.

Finally, suppose that Case (3) pertains and 
$m-km_i < 0$.  Let $q = \lfloor \frac{m-km_i}{m_i} \rfloor + 1$ (the significance of this moment 
-- ``the death of the \pep " -- was explained at the beginning of this subsection).  By the 
definition of a \pep \ (Definition \ref{pep}), the $q$-step \nf \ of $\rho_0$ admits a hard 
splitting as
\[      E_i \overline\tau^{m-qm_i} \overline\nu \gamma_{\sigma_1}^{q,m} \odot 
\sigma_{1,q}^\bullet .  \]
By Lemma \ref{pepSplitting}, the immediate future of  $E_i \overline\tau^{m-qm_i}
 \overline\nu   \gamma_{\sigma_1}^{q,m}$ admits a hard splitting as $E_i \odot \xi$.  
 Since $\gamma_{\sigma_1}^{r,m}$ is a single edge, we have a bound of $2L_1$ on the 
 length of $\xi$.  Any \nf \ of $E_i \odot \xi$ now admits a hard splitting into edge paths,
  each of which is either a \gep , a \pep \ or of length at most $V(2L_1,f)$, by induction on 
  weight and Theorem \ref{DecompTheorem}.
\end{proof}

We highlight one consequence of Proposition \ref{nfgeptopep}:

\begin{corollary}
Suppose that $\rho = E_i \overline{\tau}^{m-km_i} 
\overline{\nu} \gamma$ is a 
\pep.  Any immediate \nf \ of $\rho$ (with all nibbling on the right) 
has one of the following two forms:
\begin{enumerate}
\item $\rho' \odot \sigma$, where $\rho'$ is a \pep\ and $\sigma$ admits a hard
splitting into \atom s; or
\item $E_i \odot \sigma$, where $\sigma$ 
admits a hard splitting into \atom s.
\end{enumerate}
In particular, this is true of $f_{\#}(\rho)$.

There are entirely analogous statements in case $\rho$ is a \pep \ where
$\overline{\rho}$ has the above form and all nibbling occurs on the left.
\end{corollary}

 \section{Proof of the Beaded Decomposition Theorem}

In this section, we finally prove Theorem \ref{MainProp}.  As noted in Section 
\ref{ColourCancellation}, this immediately implies the Beaded Decomposition 
Theorem.

\begin{proof}[Proof (Theorem \ref{MainProp}).]
Take $r = k_1$, the constant from Lemma \ref{Power1}, and $J = V(2L_1,f)$, where $V$ is the constant from Theorem \ref{DecompTheorem}, and $L_1$ is the maximum length of $f_{\#}^{k_1}(E)$ for any edge $E \in G$.

Suppose that $\rho$ is an $r$-monochromatic path.  Then $\rho$ is a $1$-monochromatic path for $f_1 = f_{\#}^{k_1}$.  By Proposition \ref{gammasingleedge}, every proto-\pep \ for $f_1$ is a \pep .

By Theorem \ref{DecompTheorem}, $\rho$ admits a hard splitting into edge paths, each of which is either the \nf \ of a \gep \ or else has length at most $V(1,f)$.  By Proposition \ref{nfgeptopep}, if we replace $f$ by $f_1$ then any \nf \ of a \gep \ admits a hard splitting into edge paths, each of which is either a \gep , a \pep \ or else has length at most $V(2L_1,f)$.   By Lemma \ref{HardSplitProperties}, the splitting of the \nf \ of a \gep \ is inherited by $\rho$.

We have shown that $\rho$ is $(J,f)$-beaded, as required.
\end{proof}

\begin{remark} \label{SameJ}
We have already remarked that, for a fixed $m$, the constant $V(m,f)$ from 
Theorem \ref{DecompTheorem} remains unchanged when $f$ is replaced by an iterate.

As in Section \ref{IterateSection}, we retain the notion of $(J,f)$-beaded with the oringal
$f$ when passing
to an iterate of $f$

Therefore, when $f$ is replaced by an iterate, Theorem \ref{MainProp} remains true
with the same constant $J$.  This remark is important in our applications, for we 
pass to iterates of $f$, and the iterate chosen will depend on $J$.
\end{remark}

\section{Refinements of the Main Theorem} \label{refinements}

The Beaded Decomposition Theorem is the main result of this paper.  In this section, we 
provide a few further refinements that will be required for future applications.

Throughout this section we suppose that $f$ has been replaced with $f_1$ from 
Lemma \ref{Power1}, whilst maintaining the conventions for definitions from 
Section \ref{IterateSection}.  When we refer to $f$ we mean this iterate $f_1$.  With this in 
mind, a {\em monochromatic path} is a $1$-monochromatic path for $f$.  Similarly, 
armed with Theorem \ref{MainProp}, we refer to $(J,f)$-beads, simply as {\em beads},
and a path which is $(J,f)$-beaded will be referred to simply as {\em beaded}.
The constant $L$ now refers to the maximum length $|f(E)|$ for edges $E \in G$ { \em 
with  the new $f$}.

In the following theorem, the {\em past} of an edge is defined with respect to an 
arbitrary choice of tightening.

\begin{theorem} \label{ColourCancelMain}
There exists a constant $D_1$, depending only on $f$, with the
following properties. Suppose $i \ge D_1$, that $\chi$ is a monochromatic path and
that $\e$ is an edge in $f^{i}_{\#}(\chi)$ of weight $r$ whose past in $\chi$ is 
also of weight $r$.  Then $\e$ is contained in an edge-path $\rho$ so that $f_{\#}^i(\chi)$ has a hard splitting immediately on either side of $\rho$ and $\rho$ is one of the following:
\begin{enumerate}
\item a Nielsen path;
\item a \gep ;
\item  a \pep ; or
\item a single edge.
\end{enumerate}
\end{theorem}

\begin{proof}
Let $\chi$ be a monochromatic path.  For any $k \ge 0$, denote $f_{\#}^k(\chi)$ by $\chi_k$.  In a sense, we prove the theorem `backwards', by fixing an edge $\e_0$ of weight $r$ in $\chi_0 = \chi$ and considering its futures in the paths $\chi_k$, $k \ge 1$.  The purpose of this proof is to find a constant $D_1$ so that if $\e$ is any edge of weight $r$ in $\chi_i$ with past $\e_0$, and if $i \ge D_1$ then we can find a path $\rho$ around $\e$ satisfying one of the conditions of the statement of the theorem.

Fix $\e_0 \in \chi_0$.  By Theorem \ref{MainProp}, there is an edge path $\pi$ containing $\e_0$ so that $\chi$ admits a hard splitting immediately on either side of $\pi$ and $\pi$ either (I) is a \gep ; (II) has length at most $J$; or (III) is a \pep .  In the light of Remark \ref{stable}, it suffices to establish the existence of a suitable $D_1$ in each case.  To consider the futures of $\e_0$ in the futures $f_{\#}^k(\chi)$ of $\chi$, it suffices to consider the futures of $\e_0$ within the (entire) futures of $\pi$.  Therefore, for $k \ge 0$, let $\pi_k = f_{\#}^k(\pi)$.  Suppose that we have chosen, for each $k$, an edge $\e_k$ in $\pi_k$ such that: (i) $\e_k$ lies in the future of $\e_0$; (ii) $\e_k$ has the same weight as $\e_0$; and (iii) $\e_k$ is in the future of $\e_{k-1}$ for all $k \ge 1$.

{\bf Case (I):} $\pi$ is a \gep .  In this case, the path $\pi_k$ is a \gep \ for all $k$, any future of $\e_0$ lies in $\pi_k$, and there is a hard splitting of $\chi_k$ immediately on either side of $\pi_k$.  Therefore, the conclusion of the theorem holds in this case with $D_1 = 1$.

\medskip

{\bf Case (II):} $|\pi| \le J$.  Denote the weight of $\pi$ by $s$.  Necessarily $s \ge r$.  By Lemma \ref{SplittingLemma} the path $\pi_{D(J)}$ admits a hard splitting into edge paths, each of which is either
\begin{enumerate}
\item   a single edge of weight $s$; \label{DJ1}
\item   an indivisible Nielsen path of weight $s$; \label{DJ2}
\item   a \gep \ of weight $s$; or \label{DJ3}
\item   a path of weight at most $s-1$. \label{DJ4}
\end{enumerate}
We consider which of these types of edge paths our chosen edge $\e_{D(J)}$ lies in.  In case (\ref{DJ1})  there is a hard splitting of $\pi_{D(J)}$ immediately on either side of the edge $\e_{D(J)}$, so for all $i \ge D(J)$ there is a hard splitting of $\pi_i$ immediately on either side of $\e_i$, since $\e_i$ and $\e_{D(J)}$ both have the same weight as $\e_0$.  For cases (\ref{DJ2}) and (\ref{DJ3}), $\e_{D(J)}$ lies in an indivisible Nielsen path or \gep \ with a hard splitting of $\pi_{D(J)}$ immediately on either side, so for all $i \ge D(J)$ any future of $\e_0$ in $\pi_i$, and in particular $\e_i$, lies in an indivisible Nielsen path of \gep \ immediately on either side of which there is a hard splitting of $\pi_i$.

Finally, suppose $\e_{D(J)}$ lies in an edge path $\tilde\rho$ with a hard splitting of $\pi_{D(J)}$ immediately on either side, and that $\tilde\rho$ is not a single edge, an indivisible Nielsen path, or a \gep \footnote{In this case necessarily $s \le r-1$}.  We need only consider the future of $\tilde\rho$.  For $k \ge 0$, let $\rho_{D(J) + k} = f_{\#}^k(\tilde\rho)$ be the future of $\tilde\rho$ in $\pi_{D(J)+k}$.  Now, $|\tilde\rho| \le J L^{D(J)}$ so by Lemma \ref{SplittingLemma} the edge path $\rho_{D(J) + D(J L^{D(J)})}$ admits a hard splitting into edges paths, each of which is either
\begin{enumerate}
\item   a single edge of weight $s-1$;
\item an indivisible Nielsen path of weight $s-1$;
\item   a \gep \ of weight $s-1$; or
\item   a path of weight at most $s-2$.
\end{enumerate}

We proceed in this manner.  If we ever fall into one of the first three cases, we are done.  Otherwise, after $s-r+1$ iterations of this argument, the fourth case describes a path of weight strictly less than $r$.  Since the weight of each $\e_i$ is $r$, it cannot lie in such a path, and one of the first three cases must hold.  Thus we have found the required bound $D_1$ in the case that $|\pi| \le J$.

\medskip

{\bf Case (III):} $\pi$ is a \pep . 

Let $\pi = E_i \overline\tau^{m-km_i} \overline\nu  \gamma_{\sigma_1}^{k,m}$ as in Definition \ref{pep}.  We consider where in the path $\pi$ the edge $\e_0$ lies. First of all, suppose that $\e_0$ is the unique copy of $E_i$.  Since $\e_0$ is parabolic, it has a unique weight $s$ future at each moment in time.  Let $q = \lfloor \frac{m-km_i}{m_i} \rfloor +1$, the moment of death. For $1 \le p \le q-1$, the edge $\e_p$ is the leftmost edge in a \pep \ and there is a hard splitting of $\pi_p$ immediately on either side of this \pep .  For $p \ge q$, Lemma \ref{pepSplitting} ensures that there is a hard splitting of $\pi_p$ immediately on either side of $\e_p$. Therefore in this case the conclusion of the theorem holds with $D_1 = 1$.

Now suppose that the edge $\e_0$ lies in one of the copies of $\overline\tau$ in $\pi$, or in the visible copy of $\overline\nu$.  Then any future of $\e_0$ lies in a copy of $\tau$ or $\nu$ respectively, which lies in a \pep \ with a hard splitting immediately on either side, until this copy of $\overline\tau$ or $\overline\nu$ is consumed by $E_i$.  Again, the conclusion of the theorem holds with $D_1 = 1$.

Finally, suppose that $\e_0$ lies in $\gamma_{\sigma_1}^{k,m}$.  For ease of notation, for the remainder of the proof $\gamma$ will denote $\gamma_{\sigma_1}^{k,m}$.  By Proposition \ref{gammasingleedge} $\gamma$ is a single edge.  Until the $q$-step \nf \ of $\pi$, any future of $\gamma$ of the same weight is either $\gamma$ or will have a splitting of $\pi$ immediately on either side. 

Since $\sigma$ is an indivisible Nielsen path, and $\gamma$ is a single edge, $\gamma$ is the leftmost edge of $\overline{\sigma}$.  Therefore $[\sigma \gamma]$ is a proper sub edge-path of $\sigma$. 

Suppose that $\sigma$ has exponential weight (this weight is $r$).  By Lemma \ref{Power1} and the above remark, $f_{\#}(\sigma \gamma)$ is $r$-legal.  Therefore, any future of $\gamma$ which has weight $r$ will have, at time $q$ and every time afterwards, a hard splitting immediately on either side.

Suppose now that $\sigma$ has parabolic weight $r$.  Since $[\sigma E]$ is a proper sub edge-path of $\sigma$, and since there is a single edge of weight $r$ in $f(E)$ and this is cancelled, it is impossible for $\gamma$ to have a future of weight $r$ after time $q$.
\end{proof}

Recall that the number of strata for the map $f : G \to G$ is $\omega$.  Recall also the definition of {\em displayed} from Definition \ref{Displayed}

\begin{lemma}
Let $\chi$ be a monochromatic path.  Then the number of displayed \pep s in $\chi$ of length more than $J$ is less than $2 \omega$.
\end{lemma}
\begin{proof}
Suppose that $\chi$ is a monochromatic path, and that $\rho$ is a subpath of $\chi$, with a hard splitting immediately on either side, such that $\rho$ is a \pep , and $|\rho| > J$.  Then, tracing through the past of $\chi$, the past of $\rho$ must have come into existence because of nibbling on one end of the past of $\chi$.  Suppose this nibbling was from the left.  Then all edges to the left of $\rho$ in $\chi$ have weight strictly less than that of $\rho$, since it must have come from a proper subpath of an indivisible Nielsen path in the \nf \ of the \gep \ which became $\rho$.  Also, any \pep \ to the left of $\rho$ must have arisen due to nibbling from the left. Therefore, there are at most $\omega$ \pep s of length more than $J$ which came about due to nibbling from the left.  The same is true for \pep s which arose through nibbling from the right.
\end{proof}

\begin{lemma} \label{Displayededge}
Let $D_1$ be the constant from Theorem \ref{ColourCancelMain}, and let $f_2 = (f_1)_{\#}^{D_1}$.  If $\rho$ is an \atom , then either $(f_2)_{\#}^{\omega}(\rho)$ is a beaded path all of whose beads are Nielsen paths and \gep s, or else there is some displayed edge $\e \subseteq (f_2)_{\#}^{\omega}(\rho)$ so that all edges in $(f_2)_{\#}^{\omega}(\rho)$ whose weight is greater than that of $\e$ lie in Nielsen paths and \gep s.
\end{lemma}
\begin{proof}
Suppose that $\rho$ is an \atom \ of weight $r$. If $H_r$ is a zero stratum and $(f_2)_{\#}(\rho)$ has weight $s$ then $H_s$ is not a zero stratum.  Thus, by going forwards one step in time if necessary, we suppose that $H_r$ is not a zero stratum, so $(f_2)_{\#}(\rho)$ has weight $r$.

By Theorem \ref{ColourCancelMain}, all edge of weight $r$ in $(f_2)_{\#}(\rho)$ are either displayed or lie in Nielsen paths or \gep s (since we are considering the entire future of an \atom , \pep s do not arise here).  If all edge of weight $r$ in $(f_2)_{\#}(\rho)$ lie in Nielsen paths or \gep s then we consider the \atom s in $(f_2)_{\#}(\rho)$ of weight less than $r$ (this hard splitting exists since $\rho$ and hence $(f_2)_{\#}(\rho)$ are monochromatic paths).  We now consider the immediate future of these \atom s in $(f_2)_{\#}^2(\rho)$, etc.  It is now clear that the statement of the lemma is true.
\end{proof}

Finally, we record an immediate consequence of the Beaded Decomposition
Theorem and Proposition \ref{nfgeptopep}:

\begin{theorem}
Suppose that $\sigma$ is a beaded path.  Any \nf \ of $\sigma$ is also beaded.
\end{theorem}

\end{document}